\documentclass[12pt,a4paper,twoside]{article}
\usepackage{amssymb,amsfonts,amsmath,amsthm,euscript,mathrsfs,dsfont,bm}
\usepackage[english]{babel}
\usepackage[latin1]{inputenc}
\usepackage{dsfont}
\usepackage{color}
\usepackage{mathrsfs}
\usepackage{geometry}
\usepackage{graphicx}
\usepackage{subfigure}
\usepackage{multirow}
\usepackage{subfigure}
\usepackage{multirow, multicol}	
\usepackage{tabularx,ltxtable, booktabs}	
\usepackage{amssymb,amsfonts,amsmath,amsthm, dsfont,mathrsfs}
\usepackage{graphicx}
\usepackage{subfigure}	
\usepackage{float}	
\usepackage[bf,small]{caption}
\begin{document}


\theoremstyle{plain}
\newtheorem{thm}{Theorem}[section]
\newtheorem{cor}{Corollary}[section]
\newtheorem{prop}{Proposition}[section]
\newtheorem{lema}{Lemma}[section]
\theoremstyle{definition}
\newtheorem{exe}{Example}[section]
\newtheorem{df}{Definition}[section]
\newtheorem{rmk}{Remark}[section]
\newcommand{\Tspace}{\rule{0pt}{3ex}}
\newcommand{\Bspace}{\rule[-1.2ex]{0pt}{0pt}}
\renewcommand{\theequation}{\thesection.\arabic{equation}}
\renewcommand{\thetable}{\thesection.\arabic{table}}
\renewcommand{\thefigure}{\thesection.\arabic{figure}}
\renewcommand{\arraystretch}{1.2}
\long\def\symbolfootnote[#1]#2{\begingroup%
\def\thefootnote{\fnsymbol{footnote}}\footnote[#1]{#2}\endgroup}
\numberwithin{table}{section}
\numberwithin{figure}{section}

\newcommand{\NorT}{\ensuremath{\mathcal{N}\mathcal{T}}}
\newcommand{\Nor}{\ensuremath{\mathcal{N}}}
\newcommand{\cov}{\mathrm{Cov}}
\newcommand{\corr}{\mathrm{Corr}}
\newcommand{\C}{\mathbb{C}}
\newcommand{\E}{\mathbb {E}}
\newcommand{\F}[1]{\mathcal{F}_{#1}}
\newcommand{\iten}[1]{\vspace{.2cm}\noindent{\bf #1}}
\newcommand{\K}{\mathcal{K}}
\renewcommand{\L}{\ell}
\newcommand{\LL}{\mathcal{B}}
\newcommand{\N}{\mathbb {N}}
\newcommand{\nn}[1]{\mbox{\boldmath{$#1$}}}
\newcommand{\noise}{$\{Z_t\}_{t \in \mathbb{Z}}{}$}
\newcommand{\prob}{\mathbb{P}}
\newcommand{\pe}{$\{X_t\}_{t \in \mathbb{Z}}{}$}
\renewcommand{\proof}{\noindent\textbf{Proof: }}
\renewcommand{\qed}{\begin{flushright}\vspace{-.5cm}$\Box$\end{flushright}}
\newcommand{\R}{\mathbb {R}}
\newcommand{\stg}{$\left\{X_{t}\right\}_{t=1}^{n}{}$}
\newcommand{\var}{\mathrm{Var}}
\newcommand{\Z}{\mathbb {Z}}
\newcommand{\mat}[1]{\mbox{\boldmath{$#1$}}}
\newcommand{\fim}{\hfill $\Box$}
\newcommand{\stk}[1]{\stackrel{#1}{\longrightarrow}}
\newcommand{\s}{\hphantom{x}}
\newcommand{\Ind}[2]{\ensuremath{\mathbb {I}_{#1}{\mbox{\footnotesize $\left(#2\right)$}}}}
\newcommand{\B}{\mathcal{B}}

\thispagestyle{empty}

\vskip2cm
{\centering
\Large{\bf A Generalization of the Ornstein-Uhlenbeck Process: Theoretical Results, Simulations and Parameter Estimation}\\
}
\vspace{1.0cm}
\centerline{\large{\bf{J. Stein$^a$, S.R.C. Lopes$^b$ \symbolfootnote[3]{Corresponding author. E-mail: silviarc.lopes@gmail.com}} and A.V. Medino$^c$}}

\vspace{0.5cm}

\centerline{$^a$Federal Institute Sul-rio-grandense}
\centerline{Sapiranga, RS, Brazil}

\vspace{0.3cm}

\centerline{$^b$Mathematics Institute}
\centerline{Federal University of Rio Grande do Sul}
\centerline{Porto Alegre, RS, Brazil}

\vspace{0.3cm}

\centerline{$^c$Mathematics Department}
\centerline{University of Bras\'ilia}
\centerline{Bras\'ilia, DF, Brazil}

\vspace{0.5cm}

\centerline{\today}

\begin{abstract}
\noindent In this work, we study the class of stochastic process that generalizes the Ornstein-Uhlenbeck processes,
hereafter called by \emph{Generalized Ornstein-Uhlenbeck Type Process} and denoted by GOU type process.
We consider them driven by the class of noise processes such as Brownian motion, symmetric $\alpha$-stable
L\'evy process, a L\'evy process, and even a Poisson process. We give necessary and sufficient conditions
under the memory kernel function for the time-stationary and the Markov properties for these processes. When the
GOU type process is driven by a L\'evy noise we prove that it is infinitely divisible showing its generating triplet.
Several examples derived from the GOU type process are illustrated showing some of their basic properties as well as
some time series realizations. These examples also present their theoretical and empirical autocorrelation or
normalized codifference functions depending on whether the process has a finite or infinite second moment.
We also present the maximum likelihood estimation
as well as the Bayesian estimation procedures for the so-called \emph{Cosine process}, a particular process in the class
of GOU type processes. For the Bayesian estimation method, we consider the power series representation of Fox's H-function
to better approximate the density function of a random variable $\alpha$-stable distributed. We consider four goodness-of-fit tests for helping to decide which \emph{Cosine process}
(driven by a Gaussian or an $\alpha$-stable noise) best fit real data sets. Two applications of GOU type model are presented: one based on the Apple company stock market price data and the other based on the cardiovascular mortality in Los Angeles County data.
\newline

\vspace{0.3cm}
\vspace{.2cm}
\end{abstract}

\section{Introduction}
\renewcommand{\theequation}{\thesection.\arabic{equation}}
\setcounter{equation}{0}

In this paper, we study a class of stochastic process, given in Definition \ref{Definition_1.1}, that generalizes the Ornstein-Uhlenbeck processes, the so-called OU process. Classically, the OU process has the form given in (\ref{OUT_PROCESS}) when the integrator $L$ is the standard Brownian motion. Generalizations of this process are already known, widely studied, and applied to many areas of research. For instance, some generalizations consider the integrator $L$ in (\ref{OUT_PROCESS}) as a L\'evy process, while others a fractional Brownian motion (see \cite{Barndorff, Geulisashvili, Mai, Shao} and references therein). All these generalizations consider the exponential function as the integrand in (\ref{OUT_PROCESS}).  Definition \ref{Definition_1.1} extends this subject allowing the integrand function to come from a more general class of functions satisfying the integro-differential equation (\ref{memory_function}).

Let us present the main object of study in this paper:

\begin{df} \label{Definition_1.1} \emph{Generalized Ornstein-Uhlenbeck Type Process}

We call a stochastic process $\displaystyle V=\left(V(t); \  t\geq 0\right)$ a \emph{Generalized Ornstein-Uhlenbeck Type Process}, hereafter denoted
as \emph{GOU type process},  if it is given by
\begin{equation}\label{GOUT_process}
V(t)=V_0\rho(t)+\int_0^t\rho(t-s)dL(s), \  \forall \ t\ge 0,
\end{equation}
where
\begin{equation}\label{memory_function}
\frac{d\rho(t)}{dt} = - \int_0^t \rho(s)\,d \mu_t(s), \mbox{ with } \rho(0)=1.
\end{equation}
In (\ref{memory_function}), $\mu_t$ is a signed measure for each $t \geq 0$, $V(0)=V_0$ is called the \emph{initial condition} and it can be either
a random variable or a deterministic real constant. The function $\rho(\cdot)$ is deterministic and it is called the \emph{memory kernel function}.

In (\ref{GOUT_process}), $\displaystyle L=\left(L(t); \  t\geq 0\right)$ is the integrator process and it can be a Brownian motion,
an $\alpha$-stable process, a L\'evy process, a fractional Brownian motion, a Poisson process, a semimartingale,  or any other
class of stochastic process such that the integral in (\ref{GOUT_process}) is well defined. The process $L$ is called the
\emph{Noise Process} and $V$ is said to be driven by $L$.
\end{df}

In Definition \ref{Definition_1.1} we propose an approach to study the following generalization of the Langevin equation
 (see Kubo, 1966 and references therein)
\begin{equation}\label{gleKubo}
\left\{\begin{array}{lll}
     \frac{dV(t)}{dt} = - \int_0^t \gamma(t-s) V(s)\, ds + \eta(t)\\
      \\
     V(0)= V_0,
     \end{array}
      \right.
\end{equation}
\noindent where $\eta(\cdot)$ is a general noise process and $\gamma(\cdot)$ is a memory function.
In our approach, we transfer the problem of treating the GLE given in \eqref{gleKubo}, which would be a difficult-to-treat
integro-differential stochastic equation, to a problem of studying the more easily treatable deterministic equation given in
expression \eqref{memory_function}. In this case, the responsibility of the memory phenomena lies with the $\mu_t$ family of signed measures.
Once equation \eqref{memory_function} is solved, the resulting $\rho(\cdot)$ function is replaced into the stochastic process
\eqref{GOUT_process}.

The existence and uniqueness of the solution of equation
\eqref{memory_function} is an open problem. The examples reported here reduced to a well-known ordinary differential equation and
Theorems 2.1-2.3 in Section 2 and 3.1-3.2 in Section 3 assume the existence of a unique solution for \eqref{memory_function} and,
therefore, the process in \eqref{GOUT_process} is well defined.

Definition \ref{Definition_1.1} is an improvement of an approach presented in \cite{Medino} to study the generalized Langevin equation,
which is well known in Statistical Mechanics. Recently this GOU type process has been studied in \cite{Stein}. The GOU type process,
defined above, has several dissimilarities from O.E. Barndorff-Nielsen's theory (see \cite{BN2011}).
Here we do not assume finiteness of the second moment for
the random variable $V(t)$, for any $t\geq 0$, since we are not considering the It\^o integration.
Although, for each $t\geq 0$, the random variable $V(t)$ is infinitely divisible, the process $\displaystyle V=\left(V(t); \  t\geq 0\right)$ is neither  time-homogeneous nor a L\'evy process. In the GOU type process, the stochastic integro-differential equation generalizes both
the noise and the $V_0$ component of the process.

In this work, the authors study a class of continuous-time processes arising from
the solution of the generalized Langevin equation showing the properties of two dependence measures: codifference and spectral covariance. Their theoretical properties as well as their empirical counterparts were presented for the mentioned process in this work. These dependence measures replace the autocovariance function when it is not well defined. In \cite{Stein} the authors also proposed the maximum likelihood
estimation procedure to estimate the parameters of the process arising from the classical Langevin equation, that is, the Ornstein-Uhlenbeck
process, and of the so-called Cosine process. A simulation study for particular processes arising from this class was also proposed, giving the generation, and the theoretical and empirical counterpart for both the codifference and the spectral covariance measures.

The next examples illustrate how to construct from Definition \ref{GOUT_process} subclasses of GOU type processes.

\begin{exe}\label{Example_1.1} \emph{Ornstein-Uhlenbeck Type Process}

If in Definition \ref{Definition_1.1}, for each $t\ge 0$, we let $\mu_t= \delta_t$ be the Dirac measure concentrated at $t$ with total mass $\lambda >0$, then (\ref{memory_function}) reduces to
\begin{equation}
\rho'(t)= -\lambda \rho(t), \mbox{ with } \rho(0)=1,
\end{equation}
which has $\rho(t)=e^{-\lambda t}, \  \forall \ t\ge 0$, as the unique solution.
\noindent The corresponding GOU type process is given by
\begin{equation}\label{OUT_PROCESS}
V(t)=V_0 e^{-\lambda t}+\int_0^t e^{-\lambda(t-s)}dL(s).
\end{equation}
\end{exe}

A class of processes as in (\ref{OUT_PROCESS}) driven by L\'evy noise has been studied in the last two decades approximately as a stochastic volatility model. In this context, such stochastic processes have been called  \emph{Ornstein-Uhlenbeck type process} or in short \emph{OU type process} (see \cite{Barndorff, Geulisashvili, Mai}  and references therein). We borrow from this branch of research of Mathematical Finances our terminology of Generalized Ornstein-Uhlenbeck Type process.

When $L=W$ is the Wiener process, (\ref{OUT_PROCESS}) becomes the classical Ornstein-Uhlenbeck process. For the statistical methods in Neuroscience, this process has been applied to study the stochastic fluctuation in the membrane potential of a neuron (see \cite{Jahn, Shinomoto}, and references therein). We also mention that the process in (\ref{OUT_PROCESS}) with $L=B^H$, a fractional Brownian motion with Hurst parameter $H$, has been used in Fluid Dynamics to model turbulent flows such as homogeneous Eulerian and Lagrangian turbulence (see \cite{Shao} and references therein).

In Section $4$ we return to this example considering it with the Poisson component in the noise (see Example $4.1$).

\begin{exe}\label{Example_1.2} \emph{GOU Type Process with Sturm, Hill, Mathieu or Airy Memory Kernel}

Consider $\mu_t = \mu$ for all $t \geq 0$ and $\mu$ is absolutely continuous with respect to the Lebesgue
measure  $\lambda$, that is, $d\, \mu_t(s) = d\, \mu(s) = f(s)\, d\lambda(s)$, for all $t,s \geq 0$, where
\begin{equation}\label{RNderivative}
\nonumber
f(s) = \frac{d\,\mu}{d\,\lambda}(s)
\end{equation}

\noindent is the Radon-Nikodym derivative of $\mu$ with respect to $\lambda$. In this case, equation (\ref{memory_function}) leads to the ordinary differential equation.
\begin{equation}\label{diffequation1}
\left\{\begin{array}{ll}
    \rho''(t) + f(t)\,\rho(t)=0\\
    \rho'(0)=0, \ \ \rho(0)=1.
      \end{array}
      \right.
\end{equation}
The unique solution $\rho(\cdot)$ of this initial value problem gives the memory kernel function in the GOU type process (\ref{GOUT_process}).
\end{exe}

Ordinary differential equations as in (\ref{diffequation1}) are well known in the context of the Sturm separation theorem \cite{Simmons}. When $f(t) = a + \phi(t)$, where $a$ is a constant and $\phi(\cdot)$ is a real periodic function, equation (\ref{diffequation1}) is known as \emph{Hill's Equation}. An important particular case of Hill's equation is the \emph{Mathieu's Equation} in which $\phi(t) = b\cos(2t)$ where $b\neq 0$ is a constant (see \cite{Hale, Magnus}). Another important case of (\ref{diffequation1}) is when $f(t) = t$ and the resulting equation is called \emph{Airy Equation} \cite{Vallee} (see Example $4.7$). The authors are not aware of any application or study involving the corresponding GOU type process (\ref{GOUT_process}) in each of the sub-cases listed in this example.

In Section $4$ we return to this example considering the case when $f(t) = a^2$, with Poisson component in the noise (see Example $4.2$),
with Gaussian (see Example $4.3$) and non-Gaussian (see Example $4.5$) noise. We also consider the case when $f(t) = 2a(1-2a t^2)$,
with Gaussian (see Example $4.4$) and non-Gaussian (see Example $4.6$) noise.

The remaining part of the paper is structured as follows: GOU type processes driven by the Brownian motion is considered in Section 2. In Theorem \ref{Theorem_2.1},  we study the Gaussian structure of this class of processes; in Theorem \ref{Theorem_2.2}, we set up  necessary and sufficient conditions under the memory kernel function that says when a GOU type process with Brownian noise is time-stationary; conditions that ensure the Markov property of such process is explored in Theorem \ref{Theorem_2.3}. In Section 3, we study GOU type processes driven by non-Gaussian noise. Theorem \ref{Theorem_3.1} deals with GOU type processes driven by symmetric $\alpha$-stable noises. Under the conditions exhibited there, we prove that for each $t \ge 0$, $V(t)$ has symmetric $\alpha$-stable distribution with scale parameter $\sigma_{V(t)}$ given in (\ref{Scale_Parameter}). Theorem \ref{Theorem_3.2} considers more general L\'evy noises with given generating triplet $(G,\beta,\tau)$. In this case, we prove that $V(t)$ is infinitely divisible and its generating triplet $(A_t,\gamma_{t,v_0},\nu_t)$ is given by the set of formulae in (\ref{triplet_of_V}). Section $4$ presents a list of seven examples, some of them are
extended cases of both Examples $1.1$ and $1.2$. For all those examples we give some realization time series and, depending on the noise, we
give the theoretical and empirical autocorrelation or normalized codifference functions. In Section $5$ we present both the maximum likelihood
and the Bayesian estimation procedures for the Cosine process, studied in Examples $4.3$ and $4.5$. Also in Section $5$ we present four goodness-of-fit tests for helping to decide which \emph{Cosine process} (driven by a Gaussian or an $\alpha$-stable noise) best fit real data sets. Section $6$ presents two applications: one based on the Apple company stock market price data and the other based on the cardiovascular mortality in Los Angeles County data. Finally, Section $7$ concludes the paper.

\section{GOU Type Processes Driven by Gaussian Noise}
\renewcommand{\theequation}{\thesection.\arabic{equation}}
\setcounter{equation}{0}

In this section we consider a GOU type process $V$ as in \eqref{GOUT_process}  driven by the standard Brownian motion $\displaystyle B=\left(B(t); \ t\geq 0\right)$. That is, $V$ is given by
\begin{equation}\label{eqSOL}
V(t)=V_0 \rho(t)+\int\limits_0^t \rho(t-s) \,dB(s), \ \ t\geq 0
\end{equation}
where the function $\rho(\cdot)$ satisfies $\displaystyle \rho'(t)= - \int\limits_0^t \rho(s)\,d \mu_t(s)$, $\rho(0)=1$ and $\mu_t$ is a signed measure on the Borel $\sigma$-field of $\R$, for each $t\ge 0$.

Theorem \ref{Theorem_2.1} sets up the Gaussian structure of \eqref{eqSOL}. In Theorems \ref{Theorem_2.2} and \ref{Theorem_2.3}, we  establish necessary and sufficient conditions under the function $\rho(\cdot)$ which ensures when the process in \eqref{eqSOL} is, respectively, time-stationary and a Markov process.

\begin{thm}\label{Theorem_2.1}
Suppose that $V_0$  is a Gaussian random variable with mean 0  and variance $\sigma^2$. In addition, assume that for all $t \ge 0$, $\sigma(V_0)$ and $\sigma(B(s); \ 0 \le s \le t)$ are independent sigma-fields. Then, the GOU type process $V$ given by \eqref{eqSOL} satisfies:

\noindent \textit{(i)} For each $t\ge 0$, $V(t)$ is a zero mean Gaussian random variable with variance given by
\vspace{-0.3cm}
\begin{equation}
\displaystyle \var\left(V(t)\right)=\sigma^2\rho^2(t) + \int_0^t\rho^2(t-s)ds.
\end{equation}
\vspace{-0.8cm}

\noindent \textit{(ii)}  The autocovariance function of $V$ is given by
\begin{equation}\label{ACF}
\gamma_V(t,t+h) = \sigma^2\rho(t)\rho(t+h) + \int_0^{t}\rho(u)\rho(u+h)du, \ \forall \ t\ge 0, \ h\ge 0.
\end{equation}
\vspace{-0.5cm}
\end{thm}
\proof \textit{(i)} It is known that the stochastic integral $\displaystyle I(t) = \int\limits_0^t \rho(t-s) \,dB(s)$ has Gaussian distribution with zero mean and variance $\displaystyle \int_0^t\rho^2(t-s)ds$, as it can be seen for instance in section 4.3 of \cite{Klebaner} or section 2.3 of \cite{Kuo}. As $\sigma(V_0)$ and $\sigma(B(s); \ 0 \le s \le t)$ are independent sigma-fields, it follows that $V_0 \rho(t)$ and $I(t)$ are independent zero mean Gaussian random variables. So, it is straightforward from \eqref{eqSOL} that $V(t)$ is a zero mean Gaussian random variable with variance
\begin{align}
\displaystyle \var(V(t)) & = \displaystyle \var(V_0\rho(t)) + \var(I(t)) = \displaystyle \sigma^2\rho^2(t) + \int_0^t\rho^2(t-s)ds. \nonumber
\end{align}

\noindent \textit{(ii)} Consider the stochastic integrals $\displaystyle I(t) = \int\limits_0^t \rho(t-s) \,dB(s)$ and \newline  $\displaystyle I(t+h) = \int\limits_0^{t+h} \rho(t+h-s) \,dB(s)$. Then
\begin{align*}\label{eqACFBrowniano}
\gamma_V(t,t+h)& = \cov(V(t),V(t+h)) = \E\left[\left(V_0\rho(t) + I(t)\right)\left(V_0\rho(t+h) + I(t+h)\right)\right]\\
& = \E[V_0^2\rho(t)\rho(t+h)] + \E\left[V_0\rho(t)I(t+h)\right] + \E\left[V_0\rho(t+h) I(t)\right] + \E\left[I(t)I(t+h)\right]\\
& = \sigma^2\rho(t)\rho(t+h) + \E\left[I(t) \left(\int_0^{t}\rho(t+h-s)dB(s)+\int_{t}^{t+h}\!\!\!\rho(t+h-s)dB(s)\right)\right]\\
& = \sigma^2\rho(t)\rho(t+h) + \E\left[I(t)\!\!\!\int_0^{t}\!\!\!\rho(t+h-s)dB(s)\right] + \E\left[I(t)\!\!\!\int_{t}^{t+h}\!\!\!\!\!\!\rho(t+h-s)dB(s)\right]\\
& = \sigma^2\rho(t)\rho(t+h) + \int_0^{t}\rho(t-s)\rho(t+h-s)ds \\
& = \sigma^2\rho(t)\rho(t+h) + \int_0^{t}\rho(u)\rho(u+h)du.
\end{align*}

In the above derivation, in the fourth equality we have used that $V_0 \rho(t)$ and $I(t+h)$ are independent zero mean random variables, as well as $V_0 \rho(t+h)$ and $I(t)$ are. Second,  in the third term of line 4, the integrals $I(t)$ and $\displaystyle \int\limits_t^{t+h} \rho(t+h-s) \,dB(s)$ are independent zero mean random variables. This is a consequence of the independent increment property of the Wiener process on non-overlapping intervals, as can be seen in \cite{Klebaner, Kuo}. For the second term in line 4, we observe that $I(t)$ and $\displaystyle \int_0^{t}\rho(t+h-s)dB(s)$ are Gaussian random variables with covariance $\displaystyle \int_0^{t}\rho(t-s)\rho(t+h-s)ds$ (see  \cite{Klebaner, Kuo}). Finally, the change of variable $u=t-s$ leads to the result.
\hfill{$\square$}

\vspace{0.2cm}

Now, in Theorem \ref{Theorem_2.2}, by taking advantage of the Gaussian structure, we study the time-stationary property of the GOU type process \eqref{eqSOL} and present necessary and sufficient conditions on function $\rho(\cdot)$ which guarantees such property.

\begin{thm}\label{Theorem_2.2} Under the hypothesis in Theorem \ref{Theorem_2.1}, the process given by \eqref{eqSOL} is time-stationary if, and only if, for some $\theta>0$,
\begin{equation}
\rho(t) = e^{-\theta t}, \ \forall t\ge 0 \mbox{ and } \sigma^2 = \frac{1}{2\theta}.
\end{equation}
\end{thm}

\proof It's a well-known result that if
\begin{equation}
V(t)=V_0 e^{-\theta t}+\int\limits_0^t e^{-\theta (t-s)} \,dB(s), \ \ t\geq 0,
\end{equation}
\noindent and $\displaystyle \sigma^2 = (2\theta)^{-1}$, then $V$ is a time-stationary process. This is the classical (stationary) Ornstein-Uhlenbeck process, see for instance \cite{Karatzas}, page 358.

Suppose that $V$ is time-stationary. Then, the autocovariance function $\displaystyle \gamma_V(t,t+h)$ of $V$ does not depend on $t$ and by \eqref{ACF} we have
$$
\displaystyle \gamma_V(t,t+h)=\sigma^2\rho(t)\rho(t+h)+\psi(t,h)=\sigma^2\frac{\partial\psi}{\partial t}(t,h)+\psi(t,h),
$$
where we define $\displaystyle \psi(t,h) \equiv \int_0^t\rho(u)\rho(u+h)du$. The function $\psi(\cdot)$ is twice differentiable
since function $\rho'(\cdot)$ satisfies expression \eqref{memory_function}. So, we can take derivative with respect to $t$ on both sides of the
last above equality, we obtain that $\psi$ satisfies the partial differential equation
\[
\sigma^2\frac{\partial^2\psi}{\partial t^2}(t,h)+\frac{\partial\psi}{\partial t}(t,h)=0,
\]
with the boundary condition $\displaystyle \frac{\partial\psi}{\partial t}(0,h)=\rho(0)\rho(0+h)=\rho(h), \ h\ge 0$. Solving this equation, we have
\begin{align*}
\sigma^2\frac{\partial^2\psi}{\partial t^2}(t,h)e^{t/\sigma^2}+\frac{\partial\psi}{\partial t}(t,h)e^{t/\sigma^2}&= 0 \Longrightarrow \frac{\partial}{\partial t}\left[ \sigma^2\frac{\partial\psi}{\partial t}(t,h)e^{t/\sigma^2} \right] =0 \Longrightarrow \\
\sigma^2\frac{\partial\psi}{\partial t}(t,h)e^{t/\sigma^2}&=g(h).
\end{align*}
Taking $t=0$ and using the boundary condition, we have
\begin{align*}
\sigma^2\frac{\partial\psi}{\partial t}(0,h) = g(h)  \Longrightarrow \sigma^2\rho(h) = g(h).
\end{align*}
Hence,
\begin{align*}
\frac{\partial\psi}{\partial t}(t,h) = \rho(h)e^{-t/\sigma^2}  \Longrightarrow \rho(t)\rho(t+h) = \rho(h)e^{-t/\sigma^2} \Longrightarrow \\
\rho(t)\rho(t+0) = \rho(0)e^{-t/\sigma^2} \Longrightarrow \rho^2(t) = e^{-t/\sigma^2} \Longrightarrow \rho(t) = e^{-t/{2\sigma^2}}.
\end{align*}
So, by denoting $\displaystyle \theta = \frac{1}{2\sigma^2}$ we have the result.
\hfill{$\square$}

\vspace{0.2cm}

The next theorem is on the Markov property of the GOU type process \eqref{eqSOL}.

\begin{thm}\label{Theorem_2.3} Consider the hypothesis in Theorem \ref{Theorem_2.1} and suppose that $\rho(t)\neq 0,  \ \forall t \ge 0$ and $\mu_0(\{0\}) \ge0$. Then, the process given by \eqref{eqSOL} is  Markovian if, and only if, for some $\theta \ge 0$,
\begin{equation}\label{eqdiscrete}
\rho(t) = e^{-\theta t}, \ \forall t\ge 0.
\end{equation}
\end{thm}

\proof If $\displaystyle \rho(t) = e^{-\theta t}$ with $\theta > 0$ , then $V$ is the classical Ornstein-Uhlenbeck process. In the case of $\theta = 0$, we have $\rho(t)= 1, \ \forall t\ge 0$ and $V$ is a non-standard Brownian motion. In both cases, $V$ is a Markov process.

To prove the reciprocal, it is enough to show that the function $\rho(\cdot)$ satisfies the following exponential Cauchy equation (see \cite{Kannappan}, Corollary 1.36),
\[
\rho(t+h)=\rho(t)\rho(h), \ \forall t, h \ge 0.
\]
Recall (see \cite{Kallenberg}, Proposition 13.7) that a Gaussian process $\displaystyle X=(X(t); \ t \ge 0)$ is a Markov process if, and only if,
\begin{equation}
\displaystyle \E\left[X(s)X(u)\right]\E\left[X^2(t)\right]=\E\left[X(s)X(t)\right]\E\left[X(t)X(u)\right],  \ \forall \  0\le s \le t \le u.
\end{equation}
So, if $V$ is a Markov process, its autocovariance function $\displaystyle \gamma_V(t,t+h)$ satisfies
\begin{equation}
\displaystyle \gamma_V(0,t+h)\gamma_V(t,t) = \gamma_V(0,t) \gamma_V(t,t+h), \  \forall \  h, t\ge 0.
\end{equation}
Using \eqref{ACF}, we have
\[
\small{
\begin{array}{rcl}
\displaystyle &\sigma^2\rho(0)\rho(t+h)\left[ \sigma^2\rho^2(t)+\displaystyle \int_0^t\rho^2(u)du \right] =
\displaystyle \sigma^2\rho(0)\rho(t)\!\!\left[ \sigma^2\rho(t)\rho(t+h)\! + \!\!\!\int_0^t \!\!\!\!\! \rho(u)\rho(u+h)du \right] \Longrightarrow\\ \\
\displaystyle &\rho(t+h)\left[ \sigma^2\rho^2(t)+\displaystyle \int_0^t\rho^2(u)du \right] =
\displaystyle \rho(t)\left[ \sigma^2\rho(t)\rho(t+h)+\int_0^t\rho(u)\rho(u+h)du \right] \Longrightarrow\\ \\
\displaystyle &\sigma^2\rho^2(t)\rho(t+h)+\rho(t+h)\displaystyle \int_0^t\rho^2(u)du =
\displaystyle \sigma^2\rho^2(t)\rho(t+h)+\rho(t)\int_0^t\rho(u)\rho(u+h)du \Longrightarrow \\ \\
\displaystyle &\rho(t+h)\displaystyle \int_0^t\rho^2(u)du =
\displaystyle\rho(t)\int_0^t\rho(u)\rho(u+h)du \Longrightarrow \frac{\rho(t+h)}{\rho(t)} = \frac{\displaystyle \int_0^t\rho(u)\rho(u+h)du}{\displaystyle\int_0^t\rho^2(u)du}.
\end{array}
}
\]
Taking derivative with respect to $t$ on both sides of the last above equality, we have
\[
\small{
\begin{array}{rcl}
\displaystyle \frac{\partial}{\partial t}\left[\frac{\rho(t+h)}{\rho(t)}\right] &=& \frac{\displaystyle \rho(t)\rho(t+h)\int_0^t\rho^2(u)du - \rho^2(t)\int_0^t\rho(u)\rho(u+h)du}{\displaystyle \left(\int_0^t\rho^2(u)du\right)^2} =0.
\end{array}
}
\]
So, the function $\displaystyle \frac{\rho(t+h)}{\rho(t)} = \psi(h)$ does not depend on $t$ and
\[
\rho(t+h) = \rho(t)\psi(h), \ \forall \ h, t \ge 0.
\]
Taking $t=0$, we have $\psi(h)=\rho(h), \ \forall \ h \ge 0$. Hence,
\[
\rho(t+h)=\rho(t)\rho(h), \ \forall t, h \ge 0.
\]
So, $\displaystyle \rho(t) = e^{c\, t}$ for some real constant $c$. As it has to satisfy $\displaystyle \rho'(t)= - \int\limits_0^t \rho(s)\,d \mu_t(s), \ \ \forall t \ge 0$ and $\rho(0)=1$, we conclude that $\displaystyle c = - \int\limits_0^0 e^{c\, s}\,d \mu_0(s)= - \mu_0(\{0\}) = -\theta$, for $\theta = \mu_0(\{0\})  \ge 0$.
\hfill{$\square$}

\section{GOU Type Processes Driven by Non-Gaussian Noise}
\renewcommand{\theequation}{\thesection.\arabic{equation}}
\setcounter{equation}{0}

Following \cite{Taqqu}, see  Definition 1.1.6 and Property 1.2.5 therein, a random variable $X$ is symmetric $\alpha$-stable with index of stability $\alpha\in(0,2]$ and scale parameter $\sigma \ge 0$ if the characteristic function of $X$ is given by
\begin{equation*}
\E[e^{i \theta X}] = \exp\{-\sigma^\alpha |\theta|^\alpha\}, \  \forall \theta \in \R.
\end{equation*}
This will be denoted by $X\sim S_{\alpha}(\sigma,0,0)$. A stochastic process process $\displaystyle L=(L(t); \ t \ge 0)$ is called a symmetric $\alpha$-stable L\'evy motion if $L(0)=0\ a.s.$, $L$ has independent increments, and $\displaystyle L(t)-L(s)\sim S_{\alpha}((t-s)^{1/\alpha},0,0)$ for any $0\le s < t < \infty$. We refer the reader to \cite{Taqqu} as a standard reference to stable distributions and stable processes. Notice that the symmetric $\alpha$-stable L\'evy motion is a L\'evy process.

Theorem \ref{Theorem_3.1} sets up the stable structure of the process given by \eqref{GOUT_process} driven by symmetric $\alpha$-stable L\'evy motion.

\begin{thm}\label{Theorem_3.1}
Suppose that the GOU type process given by \eqref{GOUT_process} is driven by a symmetric $\alpha$-stable L\'evy motion $L$ with stability index
$1<\alpha<2$. Suppose that the sigma-fields $\sigma(V_0)$ and $\sigma(L(s); \ 0 \le s \le t)$ are independent, for each $t\ge 0$. If  $\ V_0\sim S_{\alpha}(\sigma_0,0,0)$, then
$$\displaystyle \E[V(t)]=0 \ \ \mbox{ and } \ \ V(t)\sim S_{\alpha}(\sigma_{V(t)},0,0),$$ where
\begin{equation}\label{Scale_Parameter}
\sigma_{V(t)}^{\alpha}=|\rho(t)|^{\alpha}\sigma_0^{\alpha}+\int_0^t|\rho(t-s)|^{\alpha}ds.
\end{equation}
\end{thm}

\proof The assumption that $1<\alpha<2$ guarantees that $\alpha$-stable distributions have finite first moment, as can be seen in property 1.2.16 in \cite{Taqqu}. Also, if $\displaystyle V_0\sim S_{\alpha}(\sigma_0,0,0)$, then property 1.2.3 in \cite{Taqqu} asserts that $\displaystyle V_0\rho(t)\sim S_{\alpha}(|\rho(t)|\sigma_0,0,0)$, in particular, $\displaystyle \E[V_0\rho(t)]=\rho(t)\E[V_0]=0$.

By denoting $\displaystyle I(t)=\int_0^t\rho(t-s)dL(s)$, we have $\displaystyle I(t)\sim S_{\alpha}(\sigma_{I(t)},0,0)$, where $\displaystyle \sigma_{I(t)}^{\alpha}=\int_0^t|\rho(t-s)|^{\alpha}ds$ (see proposition 3.4.1 in \cite{Taqqu}). So, $\E[I(t)]=0$ and $$\E[V(t)]=\E[V_0\rho(t)]+\E[I(t)]=0.$$

The independence between the sigma-fields $\sigma(V_0)$ and $\sigma(L(s); \ 0 \le s \le t)$ for each $t\ge 0$ ensures the independence between the random variables $V_0$ and $I(t)$. So, using property 1.2.1 in \cite{Taqqu}, we have
$$
V(t)=V_0\rho(t)+I(t)\sim S_{\alpha}(\sigma_{V(t)},0,0),
$$
where $\displaystyle \sigma_{V(t)}^{\alpha}=\sigma_{V_0\rho(t)}^{\alpha} + \sigma_{I(t)}^{\alpha}=|\rho(t)|^{\alpha}\sigma_0^{\alpha}+\int_0^t|\rho(t-s)|^{\alpha}ds$. \qed

 In the next theorem, we consider a GOU type process $V$, as in \eqref{GOUT_process}, driven by a L\'evy process $L$ with generating triplet $(G,\beta, \tau)$.  We recall that this means that for all $t\ge 0$, the characteristic function of $L(t)$ is given by
\begin{equation*}
\E[e^{i z L(t)}] = e^{t\psi(z)}, \  \forall z \in \R,
\end{equation*}
\noindent where
\begin{equation}\label{eqcarac}
\psi(z)=-\frac{1}{2}z^2G+i\beta z + \int_{\R}\left[e^{izy}-1-izy\,\Ind D{y}\right]\tau(dy)
\end{equation}
is the \emph{characteristic exponent} of the L\'evy process. Expression \eqref{eqcarac} is the L\'evy-Khintchine representation formula for the infinitely divisible distribution of $L(1)$. In the given generating triplet, $G\ge 0$, $\beta \in \R$, and $\tau$ is a L\'evy measure on the Borel $\sigma$-field of $\R$. Such measures are characterized by
\begin{equation}
\int_{\R} \mbox{min}\{x^2, 1 \}\tau(dx)=\int_{\R} (x^2\wedge 1) \tau(dx) < \infty \ \mbox{ and } \ \ \tau(\{0\})=0.
\end{equation}
Finally, in \eqref{eqcarac} and in the sequel, $\Ind D{\cdot}$ denotes for the indicator function of the set $D=\{x\in \R: |x|\leq 1\}$. For further details on L\'evy processes, we refer the reader to \cite{Applebaum, Sato}.

\begin{thm}\label{Theorem_3.2} Let the GOU type process in \eqref{GOUT_process} be driven by a L\'evy noise $L$ generated by the triplet $(G,\beta,\tau)$. Suppose that $\tau(\R)<\infty$, $V_0=v_0$ is deterministic and the function $\rho(\cdot)$ is continuous such that $\rho(t)\neq 0$. Then, for each $t \ge 0$, $V(t)$ is infinitely divisible and its generating triplet $(A_t,\gamma_{t,v_0}, \nu_t)$ is given by
\begin{align}\label{triplet_of_V}
\nonumber
A_t &= G\int_0^t\rho^2(t-s) ds, \\
\nonumber
\gamma_{t,v_0} &=\rho(t)v_0+\beta\int_0^t \rho(t-s)ds+\int_{\R}\tau(dy)\int_0^t y\,\rho(t-s)\left[\Ind D{y\,\rho(t-s)}-\Ind D {y}\right]ds\\
\nu_t(B) &=\int_{\R}\tau(dy)\int_0^t \Ind B {y\,\rho(t-s)}ds,\quad B\in\B(\R).
\end{align}
\end{thm}

\proof We shall consider the stochastic integration as in \cite{Applebaum}. The
following equality is true for any continuous $g(\cdot)$ function (see formula (14.71),
page 533 of \cite{Pascucci}).
\begin{equation}\label{charac_funct_Integral}
\E\left[\exp\left(i z \int_s^t g(u)dL(u)\right)\right] = \exp\left[\int_s^t\psi(z\,g(u))du\right].
\end{equation}
\noindent Taking $g(u)=\rho(u-s)$. for all $s \geq 0$, in \eqref{charac_funct_Integral}, we can write the characteristic function of $V(t)$ in (\ref{GOUT_process}) as

\begin{equation}\label{characfunctionV}
\E[e^{i z V(t)}]=\exp\left[iz\rho(t)v_0+\int_0^t \psi(z\,\rho(t-s))ds\right], \ \  z\in \R.
\end{equation}

For the term with the integral in \eqref{characfunctionV} and by using \eqref{eqcarac}, we have
\begin{align*}
\int_0^t \psi(z\,\rho(t-s))ds &= \int_0^t \left[-\frac{1}{2}(z\,\rho(t-s))^2G \,+ \right.\\
& +\left.\int_{\R}\left[e^{iz\,\rho(t-s)y}-1-iz\,\rho(t-s)y\,\Ind D{y}\right]\tau(dy)+i\beta z\,\rho(t-s) \right]ds\\
&=-\frac{1}{2}\int_0^t z^2G\rho^2(t-s)ds+ i\int_0^t\beta z\,\rho(t-s) ds +\tilde{I},
\end{align*}
where
\begin{equation*}
\tilde{I}=\int_0^t\int_{\R}\left[e^{iz\,\rho(t-s)y}-1-iz\,\rho(t-s)y\,\Ind D{y}\right]\tau(dy)ds.
\end{equation*}

We shall prove that $\tilde{I}=I$, where
\begin{equation}\label{eqI}
I=\int_{\R}\!\left(e^{izx}-1-izx\Ind D{x}\right)\nu_t(dx)+i\! \!\int_0^t\!ds\!\int_{\R}zy\,\rho(t-s)(\Ind D{y\,\rho(t-s)}-\Ind D{y})\tau(dy)
\end{equation}
and
\begin{equation}\label{eqnudx}
\nu_t(dx)=\left(\int_{\R}\int_0^t\Ind {\{y\,\rho(t-s)\}}{x}ds\tau(dy)\right)dx.
\end{equation}

It is enough to show that
\begin{align}\label{eqdifrho}
&\int_0^t ds\int_{\R}(e^{iz\,\rho(t-s)y}-1)\tau(dy)=\nonumber\\
&=\int_{\R}\left[e^{izx}-1-izx\,\Ind D{x}\right]\nu_t(dx)+i\int_0^t ds\int_{\R}zy\,\rho(t-s)\,\Ind D{y\rho(t-s)}\tau(dy).
\end{align}

Indeed, if \eqref{eqdifrho} is true, then subtracting $\displaystyle i\!\!\int_0^t ds\int_{\R}zy\,\rho(t-s)\,\Ind D{y}\tau(dy)$ from both sides of this equality, we have
\begin{align*}
\tilde I &= \int_0^t ds\int_{\R}(e^{iz\,\rho(t-s)y}-1)\tau(dy)-i\!\!\int_0^t ds\int_{\R}zy\,\rho(t-s)\,\Ind D{y}\tau(dy)=\\
&=\int_{\R}\left[e^{izx}-1-izx\,\Ind D{x}\right]\nu_t(dx)+i\!\!\int_0^t ds\int_{\R}zy\,\rho(t-s)\,\Ind D{y\,\rho(t-s)}\tau(dy)-\\
&\quad -i\!\!\int_0^t ds\int_{\R}zy\,\rho(t-s)\,\Ind D{y}\tau(dy)\Longleftrightarrow\tilde{I}=I.
\end{align*}
Given $B\in\B(\R)$, the $\sigma$-field of Borel on $\R$, using the Fubini theorem, we have
\begin{align}\label{eqnu}
\nu_t(B)&=\int_{\R}\Ind{B}{x}\nu_t(dx) =\int_{\R}\Ind{B}{x}\int_{\R}\left(\int_0^t \Ind {\{y\,\rho(t-s)\}}{x}ds\right)\tau(dy)dx\nonumber\\
&=\int_{\R}\int_{\R}\left(\int_0^t\Ind{B}{x}\Ind {\{y\,\rho(t-s)\}}{x}ds\right)\tau(dy)dx = \int_{\R}\int_0^t\Ind{B}{y\,\rho(t-s)}ds\,\tau(dy)\nonumber\\
&=\int_{\R}\tau(dy)\int_0^t\Ind{B}{y\,\rho(t-s)}ds.
\end{align}
Note that the right-hand side of equality \eqref{eqdifrho} can be written as
\begin{align*}
&\int_{\R}\left[e^{izx}-1-izx\,\Ind D{x}\right]\nu_t(dx)+i\int_0^t ds\int_{\R}zy\,\rho(t-s)\,\Ind D{y\,\rho(t-s)}\tau(dy)=A-iB+iC\\
\end{align*}
where
\begin{align*}
A & \!= \!\!\! \int_{\R} \!\!\! \left[e^{izx}-1\right] \!\! \nu_t(dx), \ \ B\!=\!\!\!\int_{\R}\!\!\!zx\,\Ind D{x}\nu_t(dx), \ \ C \! = \!\!\!\int_0^t \!\!\! ds \!\!\! \int_{\R} \!\!\! zy\,\rho(t-s)\,\Ind D{y\,\rho(t-s)} \tau(dy).
\end{align*}
Using \eqref{eqnudx} and the Fubini's theorem in $A$ and $B$, we obtain
\begin{align*}
A& \!=\!\!\! \int_{\R} \!\!\! \left[e^{izx}-1\right] \!\!\! \int_{\R} \!\!\! \left(\int_0^t \!\!\! \Ind {\{y\,\rho(t-s)\}}{x}ds\right) \! \tau(dy)dx \!= \!\!\!\int_{\R}\int_{\R} \!\!\! \left(\int_0^t \!\!\! \left[e^{izx}-1\right]\Ind {\{y\,\rho(t-s)\}}{x}ds\right) \! \tau(dy)dx\\
&=\int_0^t\int_{\R}\left(\int_{\R}\left[e^{izx}-1\right]\Ind {\{y\,\rho(t-s)\}}{x} dx\right)\tau(dy)ds =\int_0^t\int_{\R}\left[e^{izy\,\rho(t-s)}-1\right]\tau(dy)ds,
\end{align*}

\noindent and
\begin{align*}
B& \! = \!\!\! \int_{\R} \!\!\! zx\, \Ind D{x} \!\!\! \int_{\R}\left(\int_0^t \!\!\! \Ind {\{y\,\rho(t-s)\}}{x}ds\right)\!\tau(dy)dx \! = \!\!\! \int_{\R}\int_{\R}\left(\int_0^t \!\!\! zx\,\Ind D{x}\,\Ind {\{y\,\rho(t-s)\}}{x}ds\right) \! \tau(dy)dx\\
&=\int_0^t\int_{\R}\left( \int_{\R} zx\,\Ind D{x}\,\Ind {\{y\,\rho(t-s)\}}{x}dx\right)\tau(dy)ds =\int_0^t\int_{\R} zy\,\rho(t-s)\,\Ind D{y\,\rho(t-s)}\tau(dy)ds.
\end{align*}
So, $B=C$, and 	we have \eqref{eqdifrho}.

To conclude, we have to prove that $\nu_t$ is a L\'evy measure, that is,  $\displaystyle \int_{\R}(x^2\wedge 1)\nu_t(dx)<\infty$. But
\begin{equation*}
\int_{\R}(x^2\wedge 1)\nu_t(dx)=\int_{|x|\leq1}x^2\nu_t(dx)+\int_{|x|>1}\nu_t(dx).
\end{equation*}
It is enough to show that $\displaystyle \int_{|x|\leq1}x^2\nu_t(dx)<\infty$ and $\displaystyle \int_{|x|>1}\nu_t(dx)<\infty$. First, note that $x=y\,\rho(t-s)$ and $|x|> 1$ if, and only if, $|y|>\frac{1}{|\rho(t-s)|}$. Now, using the Fubini's theorem, we have
\begin{align*}
&\int_{|x|>1} \!\!\! \nu_t(dx) \! = \!\!\! \int_{|x|>1}\int_{\R}\left(\int_0^t \!\!\! \Ind {\{y\,\rho(t-s)\}}{x}ds\right) \!\tau(dy)dx \! = \!\!\!\int_0^t\int_{\R}\left( \int_{|x|>1} \!\!\! \Ind {\{y\,\rho(t-s)\}}{x}dx\right) \!\tau(dy)ds\\
&= \!\!\! \int_0^t\int_{\R}\mathbb {I}_{\{|y|>\frac{1}{|\rho(t-s)|}\}}(x)\tau(dy)ds \! = \!\!\! \int_{\R}\int_0^t\mathbb {I}_{\{|y|>\frac{1}{|\rho(t-s)|}\}}(x)ds\,\tau(dy)\\
& \leq t\int_{\R} \mathbb {I}_{\{|y|>\frac{1}{|\rho(t-s_1)|}\}}(x)\tau(dy)< \infty,
\end{align*}
where $s_1\in [0,t]$ is such that $\displaystyle \mathbb {I}_{\{|y|>\frac{1}{|\rho(t-s)|}\}}(\cdot) \le \mathbb {I}_{\{|y|>\frac{1}{|\rho(t-s_1)|}\}}(\cdot)$, for all $s\in [0,t]$.

Finally and similarly
\begin{align*}
&\int_{|x|\leq1}x^2\nu_t(dx)=\int_{|x|\leq1}x^2\int_{\R}\left(\int_0^t \Ind {\{y\,\rho(t-s)\}}{x}ds\right)\tau(dy)dx\\
&=\int_{|x|\leq1}\int_{\R}\left(\int_0^tx^2 \Ind {\{y\,\rho(t-s)\}}{x}ds\right)\tau(dy)dx =\int_0^t\int_{\R}\left(\int_{|x|\leq1}x^2 \Ind {\{y\,\rho(t-s)\}}{x}dx\right)\tau(dy)ds\\
&=\int_0^t\int_{\R}\mathbb {I}_{\{|y|\leq \frac{1}{|\rho(t-s)|}\}}(x)\rho^2(t-s)|y|^2\tau(dy)ds =\int_{\R}|y|^2\int_0^t\mathbb {I}_{\{|y|\leq \frac{1}{|\rho(t-s)|}\}}(x)\rho^2(t-s)ds\,\tau(dy)\\
&\leq t\int_{\R}|y|^2 \max_{0\leq s\leq t}\{\mathbb {I}_{\left\{|y|\leq \frac{1}{|\rho(t-s)|}\right\}}(x)\rho^2(t-s)\}\tau(dy)\\
& =t\rho^2(t-s_2)\int_{|y|\leq \frac{1}{|\rho(t-s_2)|}}|y|^2 \tau(dy)<\infty.
\end{align*}

We conclude that \eqref{characfunctionV} takes the form
\begin{align*}
&\varphi_{V(t)}(z)=\E[e^{i z V(t)}]=\exp\left\{-\frac{1}{2}z^2\int_0^t G\rho^2(t-s)ds +iz\left[\rho(t)v_0+\beta\int_0^t \rho(t-s)ds\right.\right.+\\
&\left.\left.+\int_0^t ds\int_{\R}y\,\rho(t-s)(\Ind D{y\,\rho(t-s)}-\Ind D{y})\tau(dy)\right]\right. + \left. \int_{\R}(e^{izx}-1-izx\Ind D{x})\nu_t(dx)\right\}.
\end{align*}
So, $V(t)$ is infinitely divisible with generating triplet given by the set of expressions in \eqref{triplet_of_V}.
\hfill{$\square$}

\section{Examples}
\renewcommand{\theequation}{\thesection.\arabic{equation}}
\setcounter{equation}{0}

This section is dedicated to the analysis of seven examples of processes derived from the GOU type process
defined in expression \eqref{GOUT_process}: the first one considers the OU type process (see Example $1.1$)
with Poisson component in the noise; the second one is the Cosine process with Poisson component in the noise
while Examples $4.3$ and $4.5$ revisited the Cosine process driven by, respectively, Gaussian and non-Gaussian noise.
Examples $4.4$ and $4.6$ considers the quadratic OU type process driven by, respectively, Gaussian and non-Gaussian
noise.  Finally, Example $4.7$ considers the Airy equation (see \cite{Vallee}) in equation \eqref{diffequation1}.
For all of them, we show how to generate and to simulate some of their basic properties.

Examples $4.3$ and $4.5$ will be revisited, respectively,  for both the classical estimation procedure, presented in Section 5.1,
and the Bayesian estimation procedure, presented in Section 5.2.

\begin{exe}\label{Example_4.1} \emph{OU Type Process (Example \ref{Example_1.1}) with Poisson Component in the Noise}

Consider the particular case of Example (\ref{Example_1.1}) with $\theta >0$. Let the L\'evy noise $L$ be generated by the triplet $(G,\beta, \lambda \delta_1)$, where $\delta_1$ is the Dirac measure concentrated at 1 and $\lambda >0$. For $V_0=x$ deterministic, let us apply  Theorem \ref{Theorem_3.2} to compute the elements of the generating triplet $(A_t,\gamma_{t,x},\nu_t)$ of $V(t)$.

Concerning the term $A_t$, easily we have
\begin{align}
A_t &= G\int_0^t e^{-2\theta(t-s)}ds = \frac{G}{2\theta}(1-e^{-2\theta t}).
\end{align}

Now, for the term  the $\gamma_{t,x}$
\begin{align}
\gamma_{t,x} &=\rho(t)v_0+\beta\int_0^t \rho(t-s)ds+\int_{\R}\tau(dy)\int_0^t y\,\rho(t-s)\left[\Ind D{y\,\rho(t-s)}-\Ind D {y}\right]ds\nonumber \\
&=xe^{-\theta t}+\beta\int_0^t e^{-\theta(t-s)}ds+\int_{\R}\delta_1(dy)\int_0^t y\,e^{-\theta(t-s)}\left[\Ind D{y\,e^{-\theta(t-s)}}-\Ind D {y}\right]ds\nonumber \\
&=xe^{-\theta t}+\frac{\beta}{\theta}(1-e^{-\theta t})+\int_0^t e^{-\theta(t-s)}\left[\Ind D{e^{-\theta(t-s)}}-\Ind D {1}\right]ds\nonumber \\
&=xe^{-\theta t}+\frac{\beta}{\theta}(1-e^{-\theta t}).\nonumber
\end{align}

Finally, for the Lebesgue measure, for each $B$ in the $\sigma$-field $\B(\R)$ of Borel on $\R$, we obtain
\begin{align}
\nu_t(B) &=\int_{\R}\tau(dy)\int_0^t \Ind B {y\,\rho(t-s)}ds = \int_{\R}\lambda \delta_1(dy)\int_0^t \Ind B {ye^{-\theta(t-s)}}ds = \nonumber\\
&= \lambda \int_0^t \Ind B {e^{-\theta(t-s)}}ds = \lambda \int_0^t \Ind B {e^{-\theta u}}du.
\end{align}
As for $\theta >0$, $0 < e^{-\theta u} \le 1$, for all $u\ge 0$, we have
\begin{align}
\nu_t(B) &= 0,\quad \forall B \subset (-\infty, 0]\cup (1, \infty),
\end{align}
and in general,
\begin{align}
\nu_t(B) &= \lambda \int_0^t \Ind B {e^{-\theta u}}du,\quad B\in\B(\R).
\end{align}
\end{exe}

\begin{figure}[h!!!]
\centering
\mbox{
\subfigure[]{\includegraphics[width=0.33\textwidth]{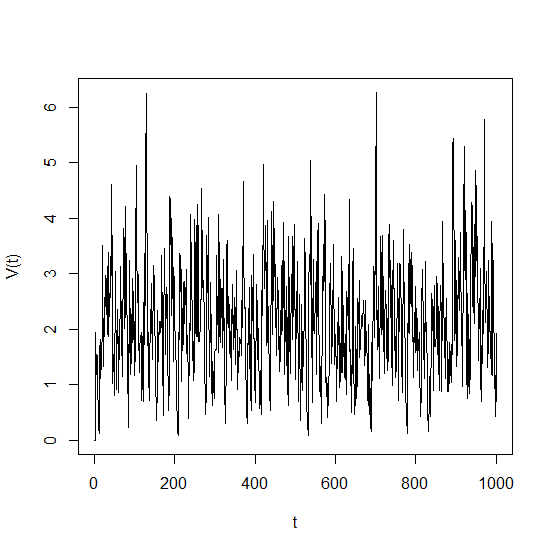}}
\subfigure[]{\includegraphics[width=0.33\textwidth]{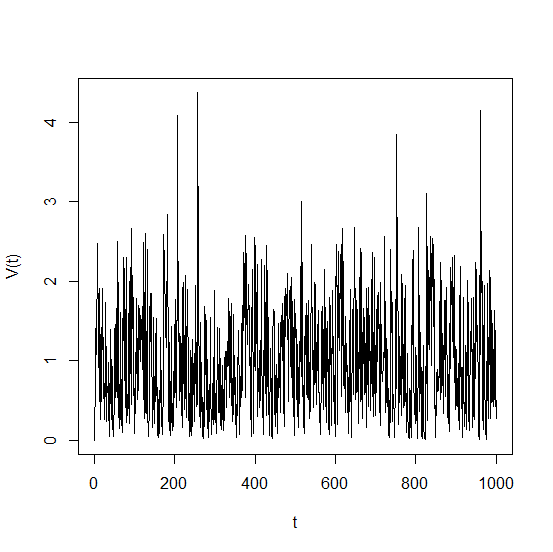}}
\subfigure[]{\includegraphics[width=0.33\textwidth]{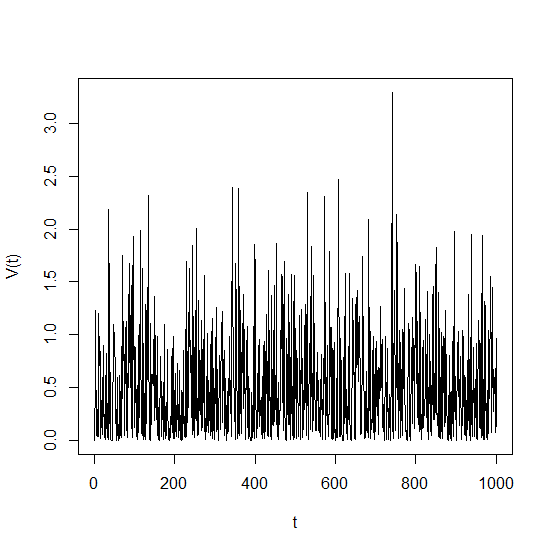}}
}
\caption{Time series derived from the OU process driven by a unit Poisson noise:
(a) $\theta=0.5$; (b) $\theta=1$; (c) $\theta=2$.}
\label{Figure_4.1}
\end{figure}

Figure \ref{Figure_4.1} shows three realization time series of the OU process in Example \ref{Example_4.1} driven by a unit Poisson noise, when the constant $\theta \in \{0.5, 1, 2\}$. In this case, the L\'evy noise $L$ has generating triplet of the form $(0, 0, 1\delta_1)$, that is, $G=\beta=0$ and $\tau = \lambda \delta_1= 1\,\delta_1$ is a Dirac measure concentrated at 1. We notice, from these three plots, that even though the noise has some jumps, the resulting process has always positive trajectories more concentrated on the abscissa axis as the value of $\theta$ increases.

\begin{exe}\label{Example_4.2} \emph{Cosine Process with Poisson Component in the Noise}

Let us recall Example \ref{Example_1.2} with the particular function $f(t) = a^2$  for all $t \geq 0$ for some constant $a>0$. Then, equation
\eqref{diffequation1} becomes
\begin{equation}
\left\{\begin{array}{ll}
    \rho''(t) + a^2\, \rho(t)=0\\
    \rho'(0)=0, \ \ \rho(0)=1.
      \end{array}
      \right.
\end{equation}
which has $\rho(t)= \cos(at), \forall \ t\ge 0$ as single solution. The respective GOU type process is
\begin{equation}\label{CosineProcess}
V(t)=V_0 \cos(at)+\int_0^t \cos(a(t-s))dL(s).
\end{equation}
It will be called \emph{Cosine OU type process}, or simply the \emph{Cosine Process}.

Let $V_0=x$ be deterministic. If this process is driven by a L\'evy noise $L$ generated by the triplet $(G,\beta,\lambda \delta_1)$, where $\delta_1$ is the Dirac measure concentrated at 1 and $\lambda >0$, to the elements of the generating triplet $(A_t,\gamma_{t,x},\nu_t)$ of $V(t)$ we have
\begin{align}
A_t &= G\int_0^t \cos^2(a(t-s))ds = G\left(\frac{t}{2}+\frac{\sin(2at)}{4a}\right).
\end{align}

For the term $\gamma_{t,x}$, we have
\begin{align}
\gamma_{t,x} &=\int_{\R}\tau(dy)\int_0^t y\,\rho(t-s)\left[\Ind D{y\,\rho(t-s)}-\Ind D {y}\right]ds+\rho(t)v_0+\beta\int_0^t \rho(t-s)ds\nonumber \\
&=\int_{\R}\delta_1(dy)\int_0^t y\,\cos(a(t-s))\left[\Ind D{y\,\cos(a(t-s))}-\Ind D {y}\right]ds +\nonumber \\
& + x \cos(at) + \beta\int_0^t \cos(a(t-s))ds = x \cos(at) + \frac{\beta \sin(at)}{a}+ \nonumber \\
&+\int_0^t \cos(a(t-s))\left[\Ind D{\cos(a(t-s))}-\Ind D {1}\right]ds =x \cos(at) + \frac{\beta \sin(at)}{a},
\end{align}
since $-1\le \cos(a(t-s)) \le 1$, for all  real numbers $st$ and $t$.

And for the Lebesgue measure, given $B\in\B(\R)$,
\begin{align}
\nu_t(B) &=\int_{\R}\tau(dy)\int_0^t \Ind B {y\,\rho(t-s)}ds = \int_{\R}\lambda \delta_1(dy)\int_0^t \Ind B {y \cos(a(t-s))}ds = \nonumber\\
&= \lambda \int_0^t \Ind B {\cos(a(t-s))}ds = \lambda \int_0^t \Ind B {\cos(au)}du = \frac{\lambda}{a} \int_0^{at} \Ind B {\cos(u)}du.
\end{align}
Since $-1\le \cos(u) \le 1$, for all $u$, we have
\begin{align}
\nu_t(B) &= 0,\quad \forall B \subset (-\infty, -1)\cup (1, \infty),
\end{align}
but in general,
\begin{align}
\nu_t(B) &= \frac{\lambda}{a} \int_0^{at} \Ind B {\cos(u)}du,\quad B\in\B(\R).
\end{align}
\begin{figure}[H]
\centering
\mbox{
\subfigure[]{\includegraphics[width=0.33\textwidth]{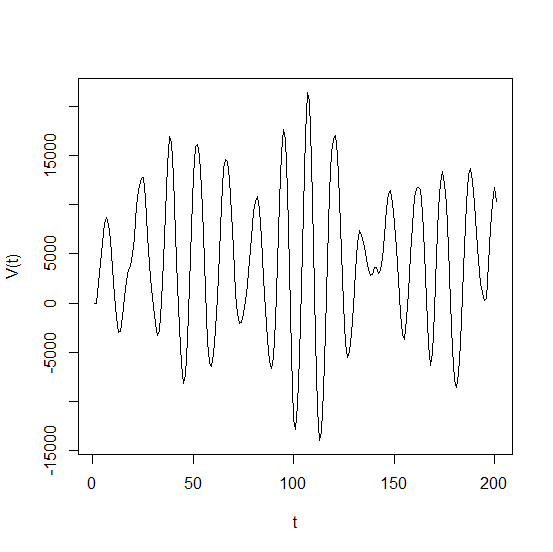}}
\subfigure[]{\includegraphics[width=0.33\textwidth]{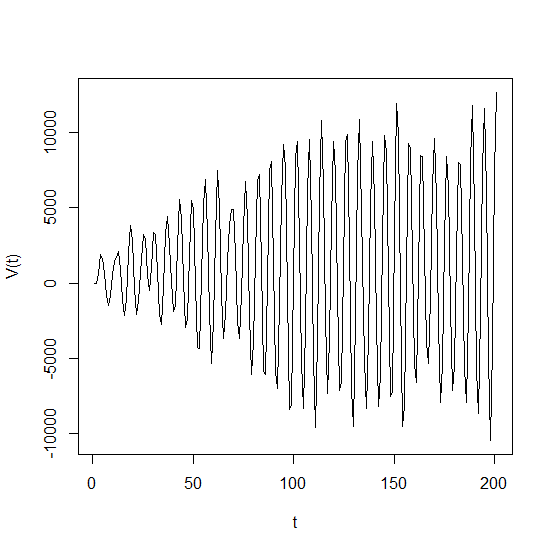}}
\subfigure[]{\includegraphics[width=0.33\textwidth]{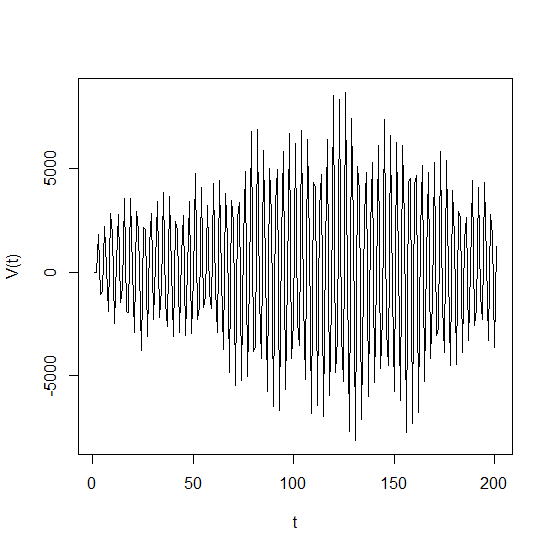}}
}
\caption{Time series derived from the Cosine process driven by a unit Poisson noise:
(a) $a=0.5$; (b) $a=1$; (c) $a=2$.}
\label{Figure_4.2}
\end{figure}

Figure \ref{Figure_4.2} shows three realization time series of the Cosine Process in Example \ref{Example_4.2} driven by a unit Poisson noise, when the constant $a \in \{0.5, 1, 2\}$. In this case, the L\'evy noise $L$ has generating triplet of the form $(0, 0, 1\delta_1)$. We notice, from these three panels, that the process maintains the periodicity feature of the cosine function. The seasonality value decreases when the value $a$ increases.

\end{exe}

\begin{exe}\label{Example_4.3} \emph{Cosine Process with Gaussian Noise}

Consider again the process \eqref{CosineProcess} but now with a Brownian noise $B$, that is
\begin{equation}\label{processCOS}
V(t)=V_0 \cos(at)+\int_0^t \cos(a(t-s))\, dB(s).
\end{equation}

\begin{figure}[H]
\begin{footnotesize}
\centering
\mbox{
   \subfigure[$a=0.5$: time series]{\includegraphics[width=0.3\textwidth]{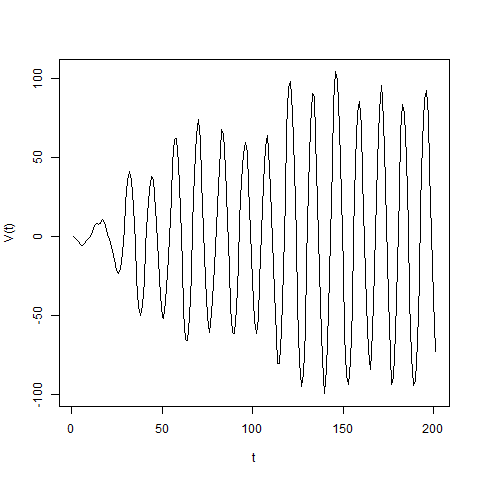}}
   \subfigure[$a=0.5$: acf (theor.)]{\includegraphics[width=0.3\textwidth]{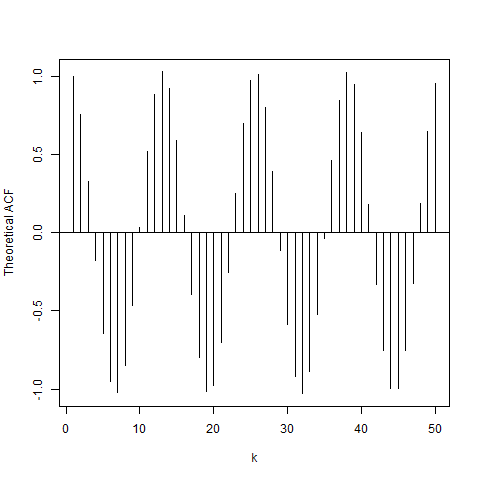}}
    \subfigure[$a=0.5$: acf (empir.)]{\includegraphics[width=0.3\textwidth]{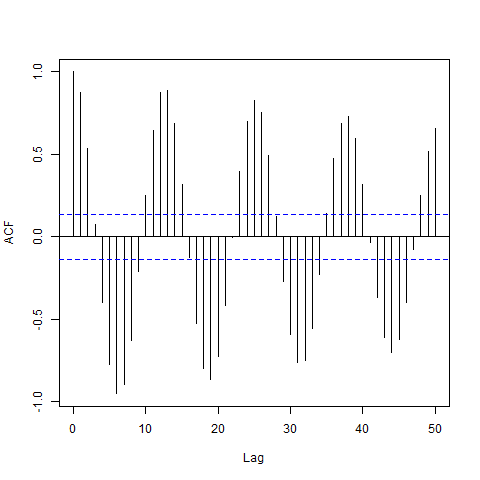}}
   }

\mbox{
   \subfigure[$a=1$: time series]{\includegraphics[width=0.3\textwidth]{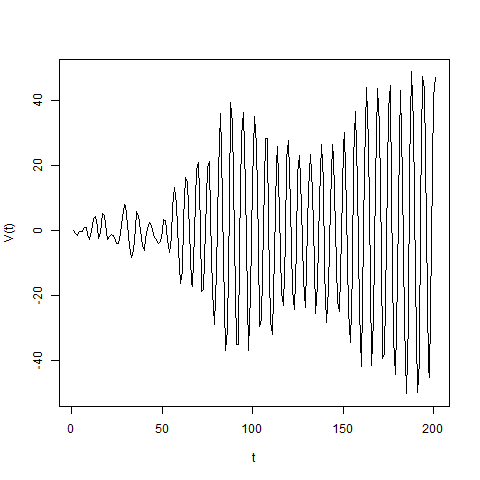}}
   \subfigure[$a=1$: acf (theor.)]{\includegraphics[width=0.3\textwidth]{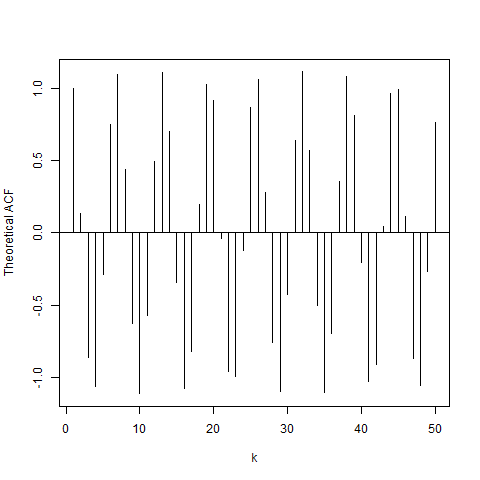}}
    \subfigure[$a=1$: acf (empir.)]{\includegraphics[width=0.3\textwidth]{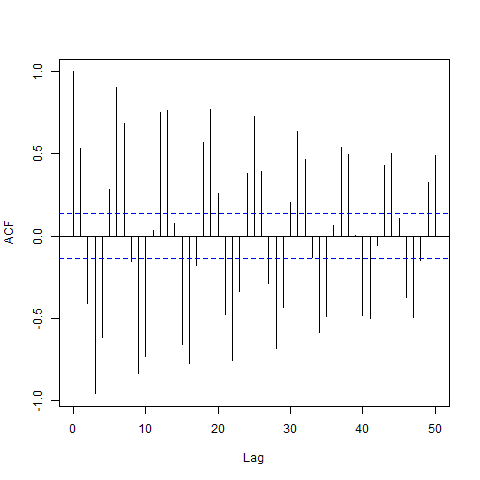}}
    }
\mbox{
   \subfigure[$a=2$: time series]{\includegraphics[width=0.3\textwidth]{serieCOSa1_alpha2_h1_n200_burn0_tfixo1_re1.png}}
   \subfigure[$a=2$: acf (theor.)]{\includegraphics[width=0.3\textwidth]{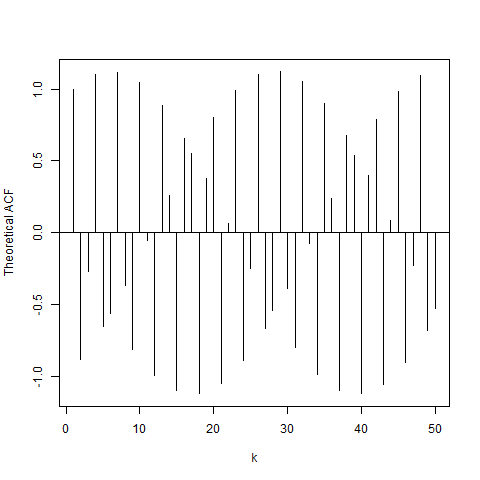}}
    \subfigure[$a=2$: acf (empir.)]{\includegraphics[width=0.3\textwidth]{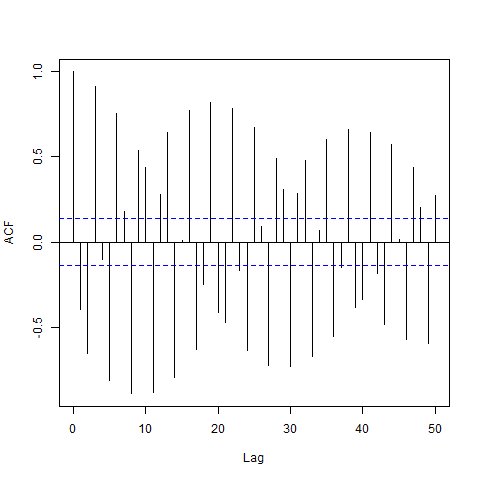}}
    }
  \caption{Time series, theoretical (theor.) and empirical (empir.) autocorrelation functions of the Cosine process, given in \eqref{processCOS}, when $a\in\{0.5,1,2\}$, $n=200$ and $h=1$.}\label{figcosACF}
  \end{footnotesize}
\end{figure}

To simulate this process we use the discrete form proposed in \cite{Stein}, given by
\begin{equation}\label{DiscreteCOS}
V((k+1)h)=2\,\cos(ah)\,V(kh)-V((k-1)h)+\varepsilon_{k,h},
\end{equation}
\noindent where $h$ is the discretization step size and $\varepsilon_{k,h}$ is a symmetric $\alpha$-stable random variable, denoted by $S_{\alpha}(\sigma_{\varepsilon, h},0,0)$, with scale parameter $\sigma_{\varepsilon, h}$ given by
\begin{equation}\label{epsilon}
\sigma_{\varepsilon, h}^{\alpha}=2\int_0^h|\cos(as)|^{\alpha}ds.
\end{equation}
As we are in the particular case of a Brownian motion noise, we have $\alpha=2$.

Figure \ref{figcosACF} shows time series, theoretical, and empirical autocorrelation functions of the Cosine process, given in \eqref{processCOS}, when $V_0\equiv 0$ and the noise is a Brownian motion. The theoretical autocorrelation function is approximated very well by its empirical counterpart. Notice that the theoretical autocorrelation function does not converge to zero, while its empirical counterpart slowly does.

\end{exe}

\begin{exe}\label{Example_4.4} \emph{Quadratic OU Type Process Driven by Gaussian Noise}

Let us consider Example \ref{Example_1.2} with $f(t)=2a(1-2at^{2})$, for $a>0$. By solving the differential equation in \eqref{diffequation1} we find $\rho(t)=e^{-at^{2}}$ and the resulting GOU type process is given by
\begin{equation}\label{processAB}
V(t)=V_0 e^{-at^{2}}+\int_0^t e^{-a(t-s)^{2}}\, dL(s).
\end{equation}
We shall call it the \emph{Quadratic OU type process}.

To simulate the process \eqref{processAB} we consider $L=B$ a Brownian noise and use the discrete form proposed in \cite{Stein}, given by

\begin{equation}\label{discreteAB}
V((k+1)h)=e^{-a(2k+1)h^2}\,V(kh)+W_{k,h},
\end{equation}
\noindent where $h$ is the discretization step size and $W_{k,h}$ is a symmetric $\alpha$-stable random variable, denoted by $S_{\alpha}(\sigma_W,0,0)$, with scale parameter $\sigma_W$ given by
\begin{equation}\label{eqsigmaAB}
\sigma_W^{\alpha}=\int_0^{kh} e^{-\alpha a((kh-s)^{2}+(2k+1)h^2)}(e^{2ash}-1)^{\alpha}\, ds+\int_{kh}^{(k+1)h} e^{-\alpha a((kh-s)^{2}-2sh+(2k+1)h^2)}\, ds.
\end{equation}
As we are in the particular case of a Brownian motion noise, we have $\alpha=2$.

Figure \ref{figQuadraticACF} shows time series, theoretical, and empirical autocorrelation functions of the quadratic process, given in \eqref{processAB}, when $V_0\equiv 0$ and the noise is a Brownian motion. Notice that, when the value $a$ increases, the theoretical function converges quickly to zero.

\begin{figure}[h!!!]
\begin{footnotesize}
\centering
\mbox{
   \subfigure[$a=0.5$: time series]{\includegraphics[width=0.3\textwidth]{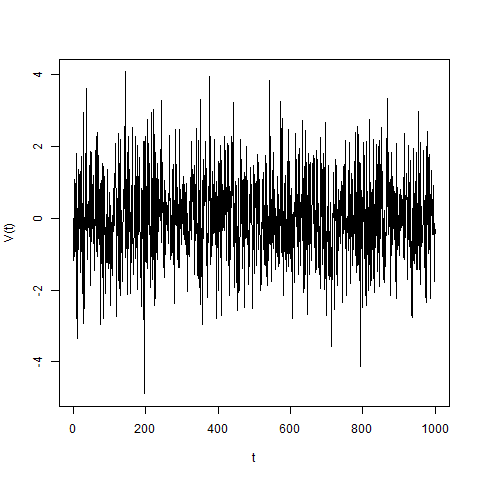}}
   \subfigure[$a=0.5$: acf (theor.)]{\includegraphics[width=0.3\textwidth]{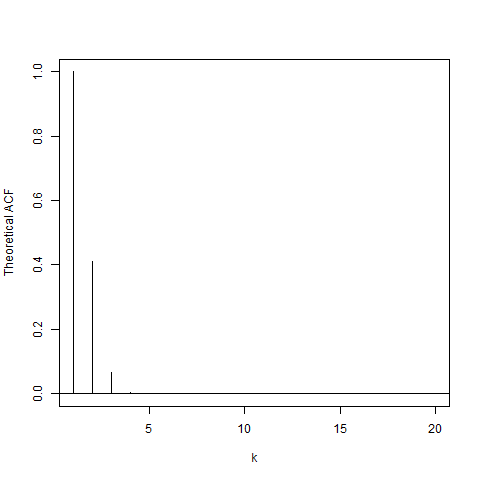}}
    \subfigure[$a=0.5$: acf (empir.)]{\includegraphics[width=0.3\textwidth]{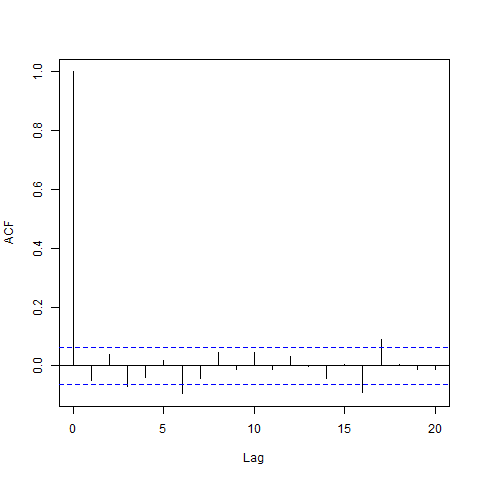}}
   }

\mbox{
   \subfigure[$a=1$: time series]{\includegraphics[width=0.3\textwidth]{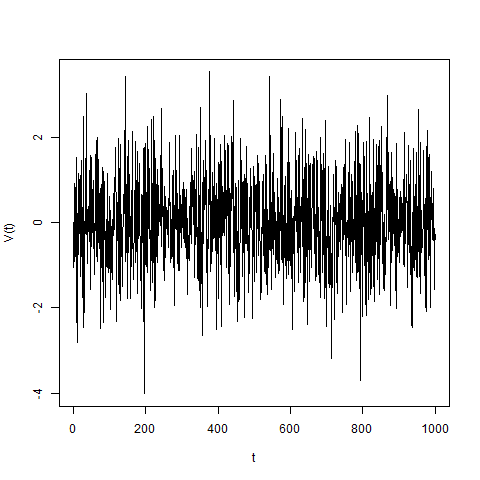}}
   \subfigure[$a=1$: acf (theor.)]{\includegraphics[width=0.3\textwidth]{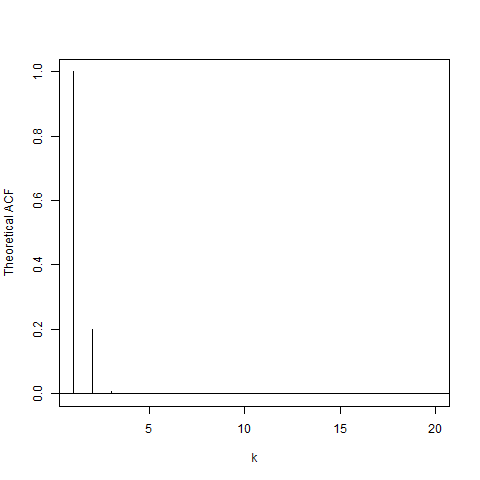}}
    \subfigure[$a=1$: acf (empir.)]{\includegraphics[width=0.3\textwidth]{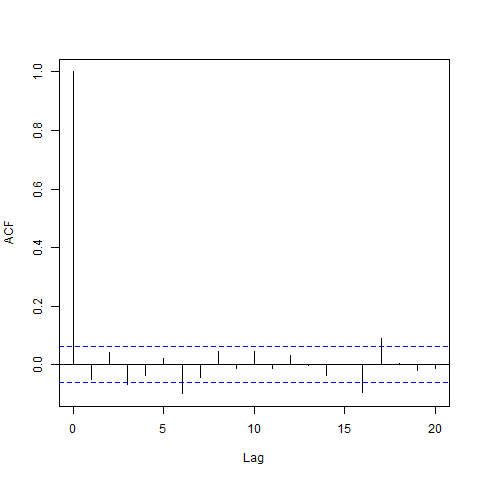}}
    }
\mbox{
   \subfigure[$a=2$: time series]{\includegraphics[width=0.3\textwidth]{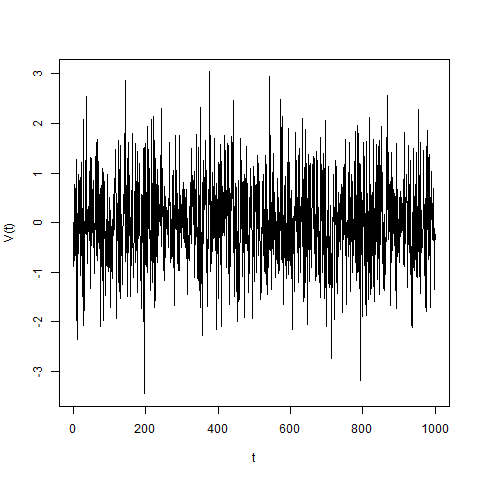}}
   \subfigure[$a=2$: acf (theor.)]{\includegraphics[width=0.3\textwidth]{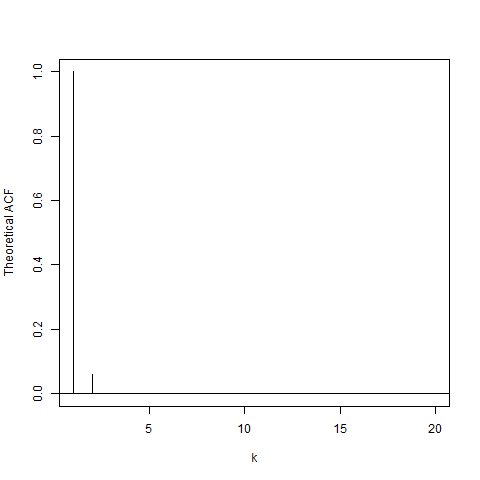}}
    \subfigure[$a=2$: acf (empir.)]{\includegraphics[width=0.3\textwidth]{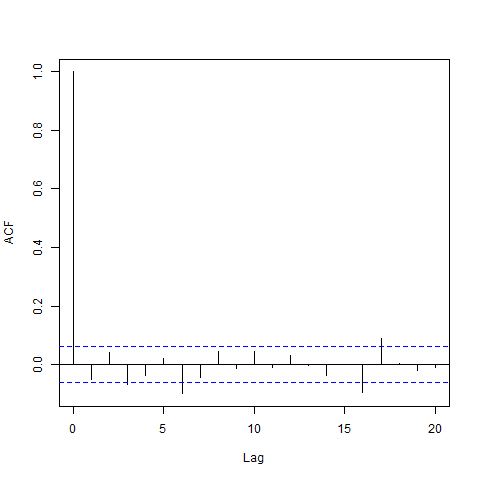}}
    }
  \caption{Time series, theoretical (theor.) and empirical (empir.) autocorrelation functions of the quadratic process, given in \eqref{processAB}, when $a\in\{0.5,1,2\}$, $n=1000$ and $h=1$.}\label{figQuadraticACF}
  \end{footnotesize}
\end{figure}

\end{exe}

\begin{exe}\label{Example_4.5} \emph{Cosine Process Driven by Non-Gaussian Noise}

Consider the Cosine process, given in \eqref{CosineProcess}, when the noise is symmetric $\alpha$-stable L\'evy motion. Since the autocovariance function is not well defined in the case of infinite second moment processes, we use the so-called codifference as a dependence measure and an estimator for it.
The codifference function is given by
\begin{align}\label{eqCOScodif1}
\tau_{V}(s;k,t)=&\ln{\left\{\mathbb{E}\left[\exp{(is(V(t+k)-V(t)))}\right]\right\}}
-\ln{\left\{\mathbb{E}\left[\exp{(is(V(t+k)))}\right]\right\}}\nonumber\\
&-\ln{\left\{\mathbb{E}\left[\exp{(-is(V(t)))}\right]\right\}},
\end{align}
\noindent where $s\in\R$, $k\geq 0$ and $t\geq 0$.

For more details, we refer the reader to \cite{Taqqu}. To avoid scale issues, we use the normalized codifference function obtained
by setting $\displaystyle \frac{\tau_{V}(s;k,t)}{\tau_{V}(s;0,t)}$, for fixed $t$.
We consider the codifference function estimator proposed in \cite{Rosadi}, for ARMA processes, given by
\begin{align}\label{eq12}
\hat{\tau}_{V}(s;k)&=\sqrt{\frac{n-k}{n}}\left[\ln{\left(\frac{1}{n-k}\sum\limits_{t=1}^{n-k} e^{is(V_{t+k}-V_t)}\right)}-\ln{\left(\frac{1}{n-k}\sum\limits_{t=1}^{n-k} e^{is V_{t+k}}\right)}\right.\nonumber\\
&\quad-\left.\ln{\left(\frac{1}{n-k}\sum\limits_{t=1}^{n-k} e^{-is V_{t}}\right)}\right],
\end{align}
\noindent where $\{V_i\}_{i=1}^n$ is a sample of size $n$ derived from the process and $k\in\{0,\cdots,n\}$. The consistency property for this estimator was derived in \cite{Stein}, for stationary symmetric $\alpha$-stable processes, with $0<\alpha\leq 2$, satisfying some conditions. The empirical normalized codifference function is obtained by setting $\displaystyle \frac{\hat{\tau}_{V}(s;k)}{\hat{\tau}_{V}(s;0)}$.

\vspace{-0.5cm}

\begin{figure}[H]
\begin{footnotesize}
\centering
\mbox{
   \subfigure[$a=0.5$: time series]{\includegraphics[width=0.3\textwidth]{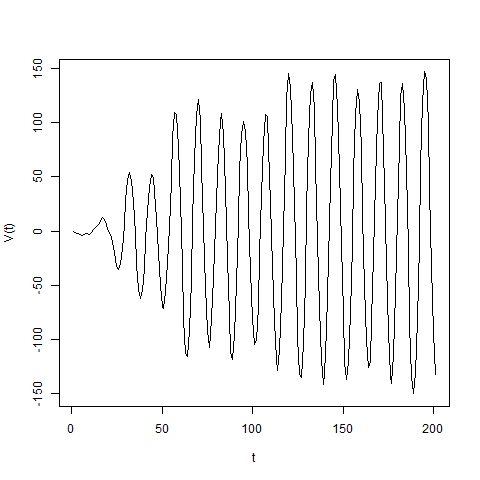}}
   \subfigure[$a=0.5$: codiff. (theor.)]{\includegraphics[width=0.3\textwidth]{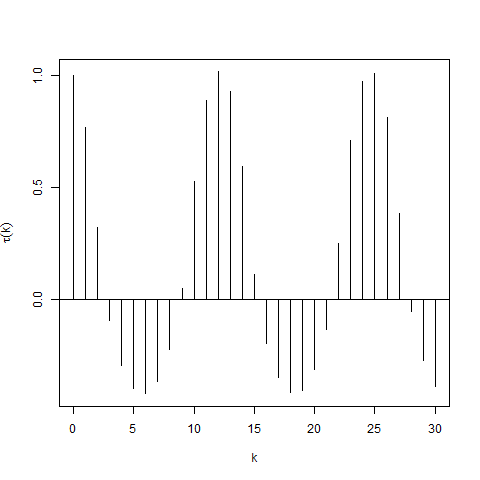}}
    \subfigure[$a=0.5$: codiff. (empir.)]{\includegraphics[width=0.3\textwidth]{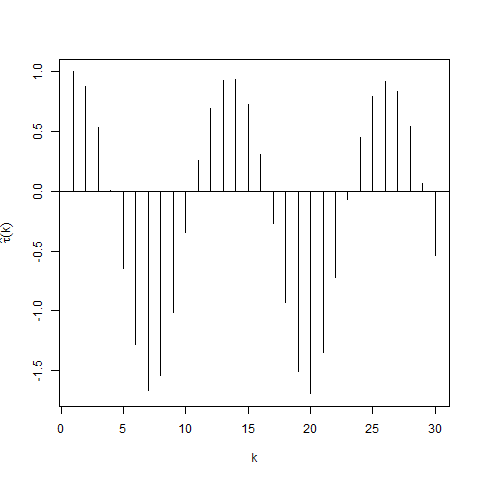}}
   }
\mbox{
   \subfigure[$a=1$: time series]{\includegraphics[width=0.3\textwidth]{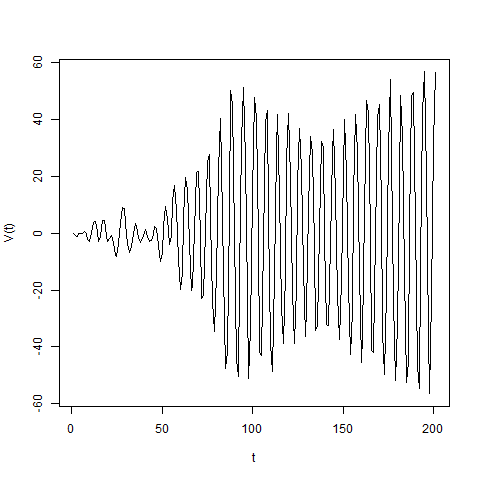}}
   \subfigure[$a=1$: codiff. (theor.)]{\includegraphics[width=0.3\textwidth]{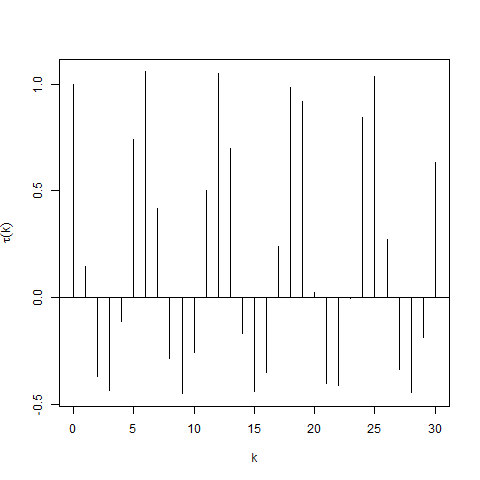}}
    \subfigure[$a=1$: codiff. (empir.)]{\includegraphics[width=0.3\textwidth]{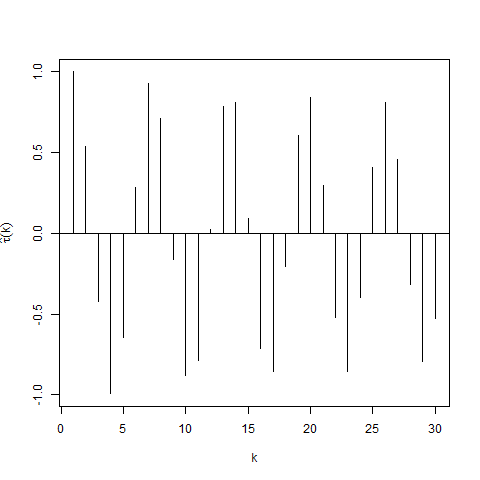}}
   }
\mbox{
   \subfigure[$a=2$: time series]{\includegraphics[width=0.3\textwidth]{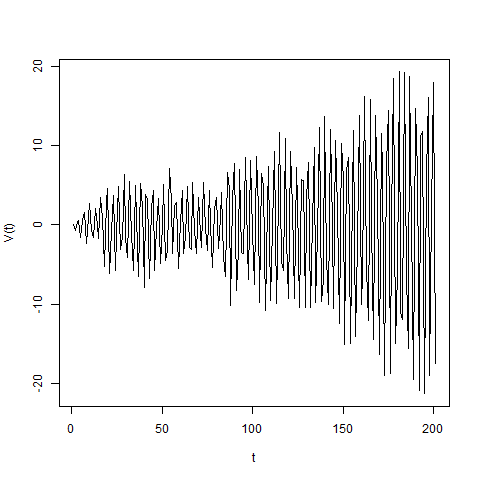}}
   \subfigure[$a=2$: codiff. (theor.)]{\includegraphics[width=0.3\textwidth]{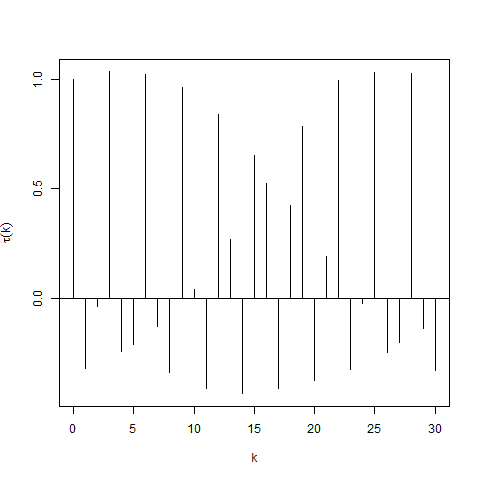}}
    \subfigure[$a=2$: codiff. (empir.)]{\includegraphics[width=0.3\textwidth]{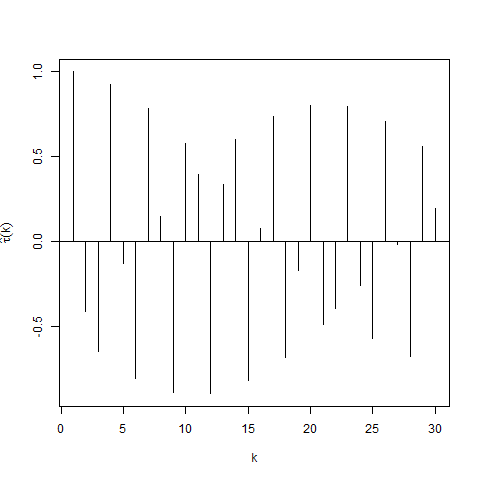}}
   }
  \caption{Time series, theoretical (theor.) and empirical (empir.) normalized codifference functions of the Cosine process,
  when the noise is a symmetric $\alpha$-stable L\'evy motion with $\alpha=1.5$, $a\in\{0.5,1,2\}$, $n=200$ and  $h=1$.}\label{figcoscodif}
  \end{footnotesize}
\end{figure}

\vspace{-0.5cm}

Figure \ref{figcoscodif} shows time series, theoretical, and empirical normalized codifference functions of the Cosine process, when $V_0\equiv 0$ and the noise is symmetric $\alpha$-stable L\'evy motion. We consider $s=0.01$ chosen to be the best choice for the $s$ value by \cite{Rosadi2}.
The codifference function of the Cosine process depends both on $k$ and $t$. Figure \ref{figcoscodif} shows the results when we fix $t=h$ in expression \eqref{eqCOScodif1}. Notice that when the value of $a$ increases, the theoretical codifference function presents large variability.

\end{exe}

\begin{exe}\label{Example_4.6} \emph{Quadratic OU Type Process Driven by Non-Gaussian Noise}

Consider the process, given in \eqref{processAB}, when the noise is symmetric $\alpha$-stable L\'evy motion. We use again
the theoretical and empirical normalized codifference functions, as described in the previous example.
Figure \ref{figABcodif} shows time series, theoretical and empirical normalized codifference functions of the quadratic process,
when $V_0\equiv 0$ and the noise is a symmetric $\alpha$-stable L\'evy motion. We consider again $s=0.01$.

\begin{figure}[h!!!]
\begin{footnotesize}
\centering
\mbox{
   \subfigure[$a=0.5$: time series]{\includegraphics[width=0.3\textwidth]{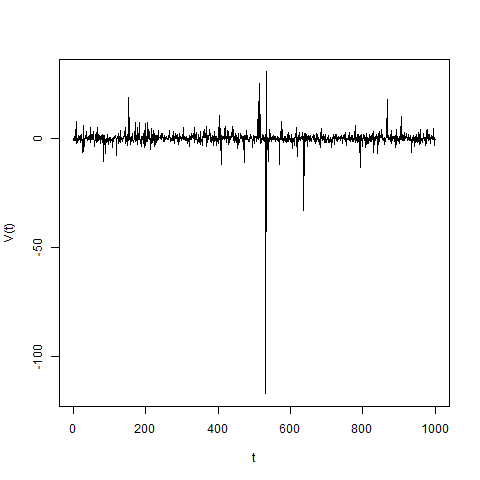}}
   \subfigure[$a=0.5$: codiff. (theor.)]{\includegraphics[width=0.3\textwidth]{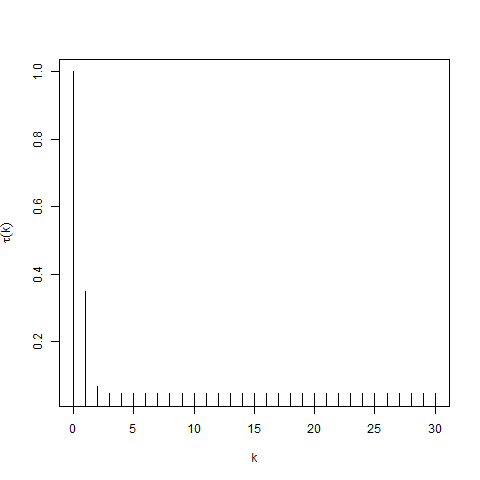}}
    \subfigure[$a=0.5$: codiff. (empir.)]{\includegraphics[width=0.3\textwidth]{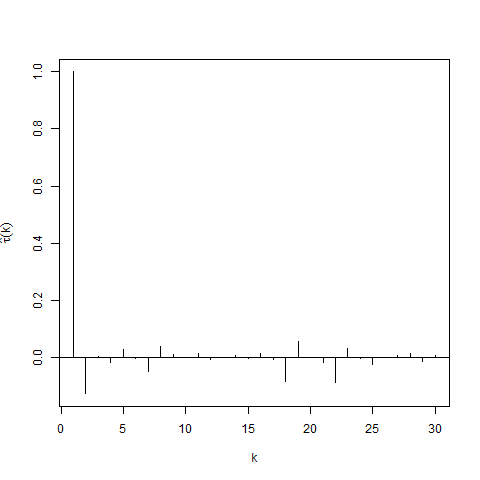}}
   }
\mbox{
   \subfigure[$a=1$: time series]{\includegraphics[width=0.3\textwidth]{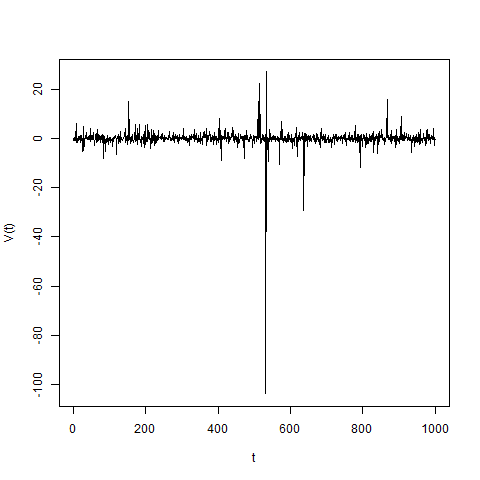}}
   \subfigure[$a=1$: codiff. (theor.)]{\includegraphics[width=0.3\textwidth]{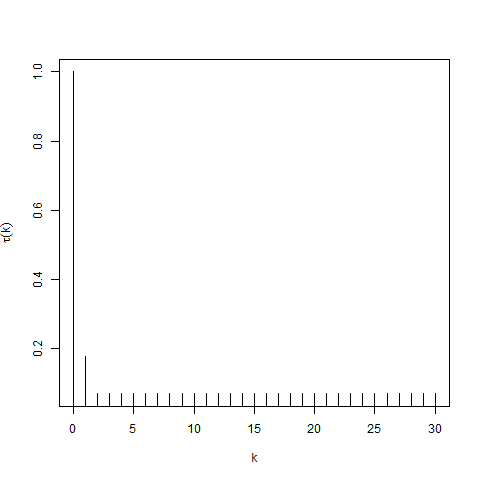}}
    \subfigure[$a=1$: codiff. (empir.)]{\includegraphics[width=0.3\textwidth]{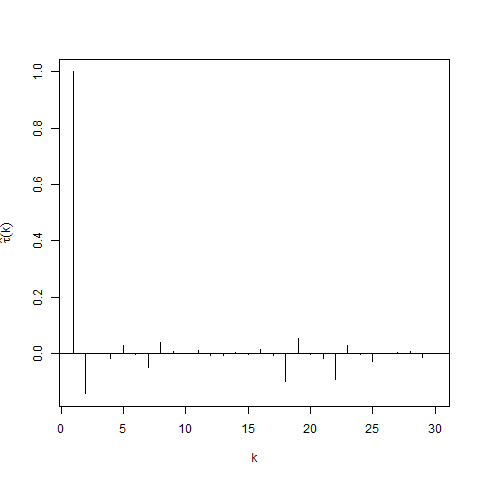}}
   }
\mbox{
   \subfigure[$a=2$: time series]{\includegraphics[width=0.3\textwidth]{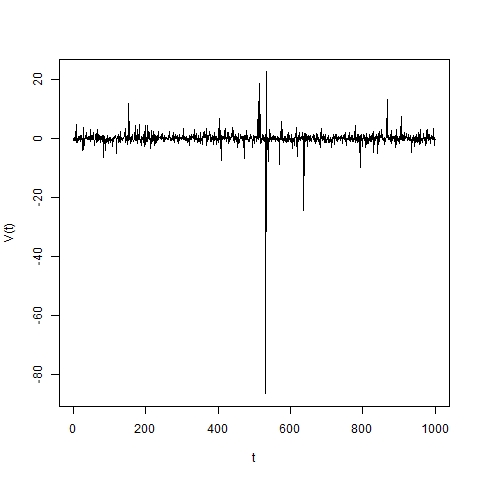}}
   \subfigure[$a=2$: codiff. (theor.)]{\includegraphics[width=0.3\textwidth]{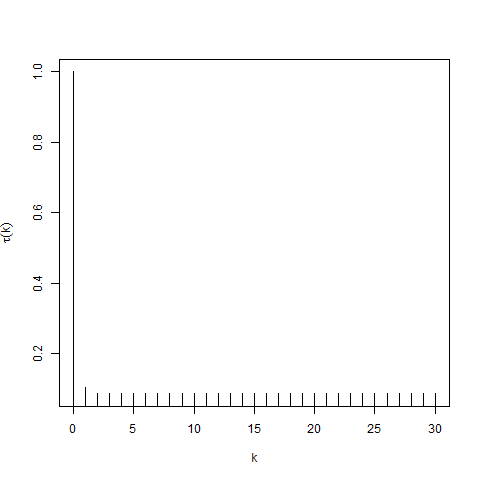}}
    \subfigure[$a=2$: codiff. (empir.)]{\includegraphics[width=0.3\textwidth]{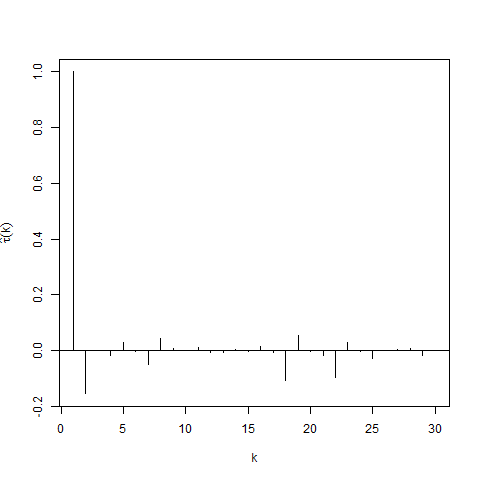}}
   }
  \caption{Time series, theoretical (theor.) and empirical (empir.) normalized codifference functions of the quadratic process,
  when the noise is a symmetric $\alpha$-stable L\'evy motion with $\alpha=1.5$, $a\in\{0.5,1,2\}$, $n=1000$ and $h=1$.}\label{figABcodif}
  \end{footnotesize}
\end{figure}
\end{exe}
\begin{exe}\label{Example_4.7} \emph{Airy Equation OU Type Process Driven by Gaussian Noise}

Let us consider Example \ref{Example_1.2} with $f(t)=t$. To solve the differential equation in \eqref{diffequation1},
we consider a power series solution of that equation given by

\begin{equation*}
\rho(t)=1+\sum\limits_{k=1}^{\infty}\frac{t^{3k}}{(2.3)(5.6).\cdots .((3k-1)(3k))}.
 \end{equation*}

To simulate the process \eqref{GOUT_process} we consider $L=B$ as a Brownian motion.
Figure \ref{Figure_example4.7} shows three realization time series of the Airy Equation OU type process driven by a Gaussian noise,
when the discretization step size is $h=0.01$ and the sample size is $n\in \{100,200,300\}$.

\begin{figure}[h!!!]
\centering
\mbox{
\subfigure[]{\includegraphics[width=0.33\textwidth]{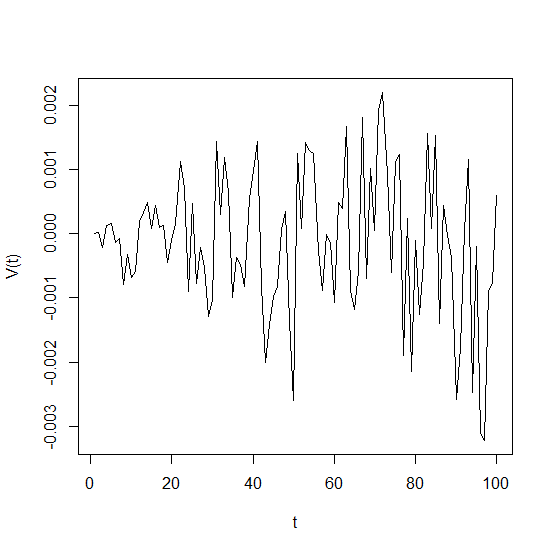}}
\subfigure[]{\includegraphics[width=0.33\textwidth]{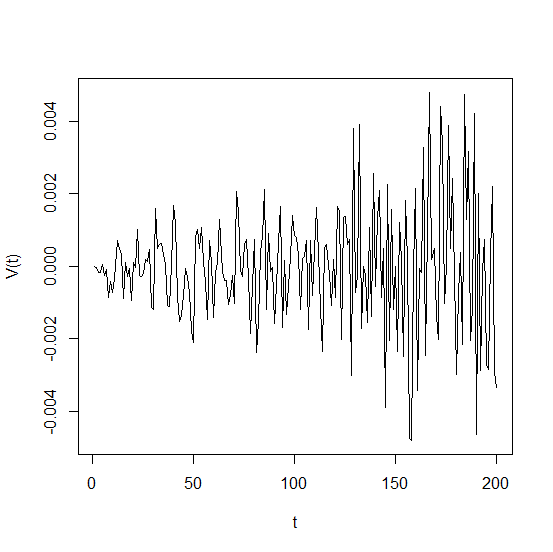}}
\subfigure[]{\includegraphics[width=0.33\textwidth]{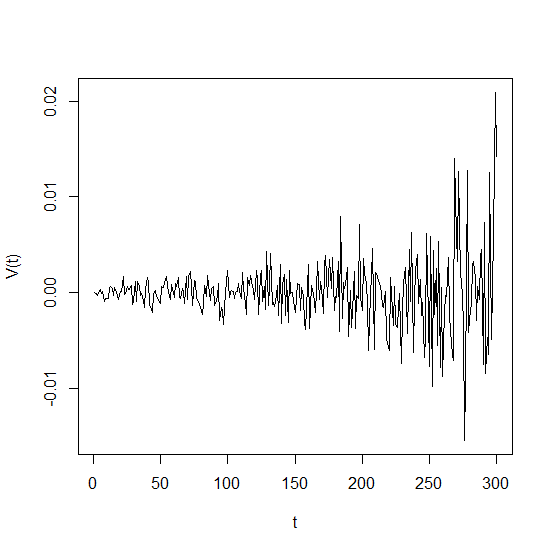}}
}
\caption{Time series derived from the Airy Equation OU type process driven by a Gaussian noise,
when the discretization step size is $h=0.01$ and the sample size is:
(a) $n=100$; (b) $n=200$; (c) $n=300$.}
\label{Figure_example4.7}
\end{figure}
\end{exe}
\section{Estimation}
\renewcommand{\theequation}{\thesection.\arabic{equation}}
\setcounter{equation}{0}

In this section, we present both the maximum likelihood (see Section 5.1) and the Bayesian (see Section 5.2) estimation procedures for the Cosine process, studied in Examples $4.3$ and $4.5$. In Section 5.3 we present four goodness-of-fit tests to be used in the Applications Section. These hypotheses tests shall be an important device to decide
which Cosine process (driven by a Gaussian or a non-Gaussian noise) best fits a real data set.

\subsection{Maximum Likelihood Estimation}

Here we present the parameter estimation, based on the maximum likelihood method, for the Cosine process driven by a Brownian motion and by symmetric $\alpha$-stable L\'evy motion, defined in Examples \ref{Example_4.3} and \ref{Example_4.5}. We point out here that we are using the maximum likelihood procedure for the processes defined in Examples 4.3 and 4.5. We know that the memory function is given by $\rho(t)=\cos(at)$, for all $t\geq 0$ and for $a$ a positive real constant. From \cite{Stein} we consider the discretization form given by the expression \eqref{DiscreteCOS}. It can be considered as a non-stationary AR$(2)$ process. Let  $\bm{\eta}=(\alpha, \sigma_{\varepsilon}, a)'$ be the parameter vector to be estimated and let $\{V_{kh}\}_{k=0}^{n-1}$ be a sample of size $n$ from the process given by the expression \eqref{processCOS}. For more details on the estimation procedure, see \cite{Stein}.

Table \ref{figMLECOS1} presents the results for the maximum likelihood estimation based on Monte Carlo simulations. With the following scenarios
\[
\bullet \ \ a\in\{1,2\}; \ \ \  \bullet  \ \ \alpha \in\{1.1,1.5,2\}; \ \ \ \bullet \ \ n=2000;  \ \ \ \bullet  \ \ h=1;  \ \ \ \bullet \ \ 500 \mbox{ replications}.
\]
From Table \ref{figMLECOS1} we observe that the maximum likelihood estimation procedure has a very good  performance. For each parameter in the vector $\bm{\eta}=(\alpha, \sigma_{\varepsilon}, a)'$, this table reports  the mean, the bias, the standard error (s.e.) and the confidence interval (CI) for its estimator. We observe that the estimation for the $a$ parameter  is very accurate. The estimation procedure for all three parameters  improves when $\alpha$ approaches to $2$, the Gaussian case. By comparing all the results, it is clear that the worst scenario occurs  when $\alpha=1.1$, since the bias is slightly greater than the biases for the other cases. There is no significant difference among the results obtained for the cases $a=1$ or $a=2$.
\begin{table}[!h]
\caption{Maximum likelihood estimation results for the Cosine process, when $a\in\{1,2\}$,
$\alpha\in\{1.1,1.5,2\}$, $n=2000$, $h=1$ and $500$ replications.}\label{figMLECOS1}
\vspace{-0.5cm}
\begin{center}
{\footnotesize
\begin{tabular}{|c|ccc|}
\hline
\hline
Statistic   & $\alpha=1.1$  & $\sigma_\varepsilon=1.5824$ & $a=1$     \\
\hline
\hline

      $\mbox{mean}$ &     1.1294 &     1.6054 &     1.0000  \\

      $\mbox{bias}$ &    -0.0294 &   -0.0230  &    -2.8e-06\\

      $\mbox{s.e.}$ &    0.0410 &     0.0679   &     0.0001  \\

      $\mbox{CI}$ &    [1.0903, 1.1663] &  [1.5014, 1.7132]  & [0.9999, 1.0001] \\
\hline
\hline
Statistic    & $\alpha=1.1$ & $\sigma_\varepsilon=1.0452$ & $a=2$       \\
\hline

      $\mbox{mean}$ &    1.1270 &       1.0586  &    1.9999 \\

      $\mbox{bias}$ &    -0.0270 &  -0.0134  &    1.6e-06\\

      $\mbox{s.e.}$ &    0.0436 &     0.0476   &     0.0001  \\

      $\mbox{CI}$ &    [1.0886, 1.1638] &  [0.9865, 1.1268]  & [1.9999, 2.0001] \\
\hline
\hline
Statistic    & $\alpha=1.5$ & $\sigma_\varepsilon=1.3450$ & $a=1$    \\
\hline

       $\mbox{mean}$ &     1.5016 &      1.3446 &    1.0000 \\

      $\mbox{bias}$ &    -0.0016 &   0.0004 &    -6.1e-07\\

      $\mbox{s.e.}$ &    0.0383 &     0.0337  &     0.0003  \\

      $\mbox{CI}$ &    [1.4345, 1.5654] &  [1.2779, 1.4066]  & [0.9994, 1.0005] \\
\hline
\hline
Statistic   & $\alpha=1.5$ & $\sigma_\varepsilon=0.9467$ & $a=2$    \\
\hline

      $\mbox{mean}$ &     1.5015 &      0.9481  &   2.0000 \\

      $\mbox{bias}$ &    -0.0015 &   -0.0014 &    -1.5e-06 \\

      $\mbox{s.e.}$ &    0.0406 &    0.0392 &     0.0003  \\

      $\mbox{CI}$ &    [1.4348, 1.5658] &  [0.8986, 0.9901]  & [1.9994, 2.0006] \\
\hline
\hline
Statistic & $\alpha=2$ & $\sigma_\varepsilon=1.2061$ & $a=1$    \\
\hline

      $\mbox{mean}$ &     1.9999 &     1.2045  &   0.9999 \\

       $\mbox{bias}$ &    0.0001 &   0.0016 &    1.0e-05\\

      $\mbox{s.e.}$ &    0.0015 &    0.0194 &     0.0007  \\

      $\mbox{CI}$ &   [1.9999,2] &  [1.1671, 1.2405]  & [0.9985, 1.0014] \\
\hline
\hline
Statistic & $\alpha=2$ & $\sigma_\varepsilon=0.9004$ & $a=2$    \\
\hline

      $\mbox{mean}$ &     2.0000 &     0.8993  &   1.9999 \\

      $\mbox{bias}$ &    0.0000 &   0.0011 &    6.5e-07\\

      $\mbox{s.e.}$ &    0.0005 &    0.0144 &     0.0007  \\

      $\mbox{CI}$ &    [1.9999, 2.0000] &  [0.8718, 0.9260]  & [1.9985, 2.0014] \\
\hline
\hline
\end{tabular}}
\end{center}
\end{table}

\subsection{Bayesian Estimation}

In this subsection we present the Bayesian methodology to estimate all three parameters $\bm{\eta}=(\alpha, \sigma_{\varepsilon}, a)'$
for the Cosine process. We point out here that we are using the Bayesian methodology for the processes defined in Examples 4.3 and 4.5,
based on Fox's H-function series representation given in this subsection.
We now consider both the Brownian motion and the symmetric $\alpha$-stable L\'evy motion.
For the Bayesian methodology we use the software \emph{JAGS} for the Gibbs sampler (see, respectively, \cite{Plummer} and \cite{Gelfand})
together with \emph{R} software, through the package \emph{R2jags}.

To obtain the posterior distribution we often deal with some complicated integrals that have no closed-form. To work around this
problem an alternative way is to consider Monte Carlo methods via Markov Chain (MCMC). In this context, we may generate
a sample from the posterior distribution in order to calculate estimates of interest for this distribution.
The two most useful MCMC methods are the Gibbs sampler and the Metropolis-Hastings algorithms. Here we consider the Gibbs sampler
to generate samples from the posterior distributions.

It is well-known that the $\alpha$-stable distribution has closed-form only for three values of $\alpha$: in the L\'evy case (when
$\alpha=0.5$), in the Cauchy case (when $\alpha=1$) and in the Gaussian case (when $\alpha=2$). To overcome the problem of the numerical
integration for others values of $\alpha$ faced when using Bayesian estimation methods we first have considered a proper bivariate probability density
of $(X,Y)$, where $X$ is $\alpha$-stable distributed and $Y$ is an auxiliary random variable such that the joint density function
has a known form. This idea was considered in the paper \cite{Buckle} where the posterior distribution can formally
be obtained through the Bayes's theorem by integrating out the unwanted random variable $Y$.
Unfortunately, through this method, the chains did not converge and the samples did present a strong correlation feature leading to bad results. We have decided to consider a power series approximation for the $\alpha$-stable density by the
so-called Fox's H-function or the H-function, involving Mellin-Barnes integrals, which is a generalization of the G-function of Meijer.
For more details on the H-function, we refer the reader to \cite{Mathai} and \cite{Rathie}.

The H-function is defined by means of a Mellin-Barnes type integral
\begin{eqnarray}\label{Hfunc}
\nonumber
H(z)& = &H_{p,q}^{m,n}\left[\begin{array}{c|ccccccc}
{\multirow{2}{*}{$z$}} & (a_{1},A_{1}), & \ldots, & (a_{n},A_{n}), & (a_{n+1},A_{n+1}), & \ldots, & (a_{p},A_{p})  \\
& (b_{1},B_{1}), & \ldots, & (b_{m},B_{m}), & (b_{m+1},B_{m+1}), & \ldots, & (b_{q},B_{q})\\
\end{array}
\right]\\
&=& \frac{1}{2\pi i}\int_{L} \frac{\displaystyle\prod_{j=1}^{m}\Gamma(b_j +B_{j}s)\displaystyle\prod_{\ell=1}^{n}\Gamma(1 - a_\ell - A_{\ell}s)}{\displaystyle\prod_{j=m+1}^{q}\Gamma(1 - b_j - B_{j}s)\displaystyle\prod_{\ell=n+1}^{p}\Gamma(a_\ell +A_{\ell}s)}z^{-s}ds,
\end{eqnarray}
\noindent where $A_\ell$ and $B_j$ are positive real constants while $a_\ell$ and $b_j$ can be real or complex-valued constants, for all $j=1, \cdots, q$
and $\ell=1, \cdots, p$. The suitable contour $L$ separates the poles of the Gamma function $\Gamma(b_j + B_{j}s)$, for $j = 1, \cdots, m$ from the poles
of the Gamma function $\Gamma(1 - a_\ell + A_{\ell}s)$, for $\ell = 1, \cdots, n$.

The H-function can be approximated by a power series (see \cite{Braaksma} and \cite{Rathie}).
For $z \neq  0$,  if $\kappa > 0$ or for $0 < |z| > D^{-1}$ if $\kappa = 0$, we have
\begin{eqnarray}\label{Hpower1}
\nonumber
H_{p,q}^{m,n}(z) &=& \displaystyle\sum_{h=1}^{m}\displaystyle\sum_{v=0}^{\infty}\frac{\displaystyle\prod_{j=1,  j \neq h}^{m}\Gamma\left(b_j +B_{j}\frac{b_h+v}{B_h}\right)}{\displaystyle\prod_{j=m+1}^{q}\Gamma\left(1 - b_j +B_{j}\frac{b_h + v}{B_h}\right)} \\
&\qquad\times&
\frac{\displaystyle\prod_{\ell=1}^{n}\Gamma\left(1 - a_\ell + A_{\ell}\frac{b_h+v}{B_h}\right)}{\displaystyle\prod_{\ell=n+1}^{p}\Gamma\left(a_\ell - A_{\ell}\frac{b_h + v}{B_h}\right)} \frac{(-1)^{v}z{\frac{b_h + v}{B_h}}}{v!B_h}.
\end{eqnarray}

For $z \neq  0$,  if $\kappa < 0$, or for $|z| > D^{-1}$, if $\kappa = 0$, we have
\begin{eqnarray}\label{Hpower2}
\nonumber
H_{p,q}^{m,n}(z) &=& \displaystyle\sum_{h=1}^{n}\displaystyle\sum_{v=0}^{\infty}\frac{\displaystyle\prod_{\ell=1,  \ell \neq h}^{n}\Gamma\left(1 - a_\ell +A_{\ell}\frac{1- a_h + v}{A_h}\right)}{\displaystyle\prod_{\ell=n+1}^{p}\Gamma\left(a_\ell +A_{\ell}\frac{1 - a_h + v}{A_h}\right)}\\
&\qquad\times&\frac{\displaystyle\prod_{j=1}^{m}\Gamma\left(b_j + B_{j}\frac{1 - a_h + v}{A_h}\right)}{\displaystyle\prod_{j=m+1}^{q}\Gamma\left(1 - b_j - B_{j}\frac{1 - a_h + v}{A_h}\right)} \frac{(-1)^{v}(1/z)^{\frac{1 - a_h + v}{A_h}}}{v!A_h},
\end{eqnarray}
\noindent where
\begin{equation*}
\kappa = \displaystyle\sum_{j=1}^{q}B_j - \displaystyle\sum_{\ell=1}^{p}A_\ell \qquad \mbox{and} \qquad D = {\displaystyle\prod_{\ell=1}^{p}A_\ell^{A_\ell}}\times{\displaystyle\prod_{j=1}^{q}B_j^{-B_j}}.
\end{equation*}

By using the power series representation of the H-function, given in expressions \eqref{Hpower1} and \eqref{Hpower2},
we can approximate the density function of a random variable $\alpha$-stable distributed.
We want to give the expression of $f_X(\cdot)$ when $X \sim S_{\alpha}(\sigma,\beta,\mu)$, considering the formula given in \cite{Schneider}.

First of all, we need to verify what is the parameterization used by \cite{Schneider}. From the expressions (2.1)-(2.2) in \cite{Schneider} considering the notation of his work, we obtain
\begin{equation} \label{Psi}
\Psi_{\alpha, \beta}(k) = -k^\alpha exp\left\{\frac{-i\pi}{2}\left(\frac{-\beta}{K(\alpha)}\right)K(\alpha) sign(k) \right\}, \ k \ge 0,
\end{equation}
where $\Psi_{\alpha, \beta}(k) = \overline{\Psi_{\alpha, \beta}(-k)}$, for $k<0$, $K(\cdot)$ defined as

\begin{equation}\label{Kalpha}
K(\alpha) = \alpha - 1 + sign(1- \alpha) = \left\{\begin{array}{ll}
    \alpha, \quad\mbox{if}\quad 0<\alpha<1
\\
    \alpha - 2, \quad\mbox{if}\quad 1<\alpha<2,
    \end{array}
\right.
\end{equation}
\noindent and $sign(\cdot)$ is the sign function.

Expression \eqref{Psi} above coincides with the parameterizarion (B) of \cite{Zolotarev} when
\[
0 = \lambda \gamma \equiv \mu; \ \ 1 = \lambda \equiv \sigma^\alpha; \ \ \beta' = -\beta / K(\alpha) \ \mbox{ and } \ \alpha \neq 1.
\]
The author \cite{Schneider} does not consider $\alpha = 1$. Besides, in his work he considers
\begin{equation} \label{Psi1}
0 < \alpha < 1, \ \ |\beta| \le \alpha \ \Longrightarrow \ |\beta ' | = |-\beta / K(\alpha)| \le \alpha / \alpha = 1
\end{equation}
and
\begin{equation} \label{Psi2}
1 < \alpha < 2, \ \ |\beta| \le 2-\alpha \ \Longrightarrow \ |\beta ' | = |-\beta / K(\alpha)| \le (2-\alpha) /(2- \alpha) = 1.
\end{equation}
Hence, from \eqref{Psi1} and \eqref{Psi2} the symmetry parameter $\beta'$ needs to be $\beta '  = -\beta / K(\alpha)$,
where $\beta$ is the symmetry parameter of \cite{Schneider}, in order to use his formula (2.17),
for $\alpha >1$, or (2.18), for $\alpha <1$.

Recall that \cite{Schneider} considers $X \sim S_{\alpha}(1,\beta,0)$ stable random variable and
$f_{\alpha,\beta}(\cdot)$ means the density function $f_X (\cdot,\alpha, 1,\beta,0)$.

We now will consider two cases: $\alpha >1$ and $\alpha< 1$.

\vspace{0.3cm}

\noindent\textbf{Case 1: $\alpha >1$}

Expression (2.17) in \cite{Schneider} gives us
\begin{equation*} \label{fX}
f_X (x;\alpha, 1,\beta,0)= \frac{1}{\pi}\sum_{n=1}^\infty \frac{\Gamma(1+n\epsilon)}{n!} \sin{(\pi n \gamma)}(-x)^{n-1},
\end{equation*}
for $\alpha >1$, where $\epsilon = 1/\alpha$ and $\gamma = (\alpha - \beta)/(2\alpha) = (\alpha + \beta'K(\alpha))/(2\alpha)$,
with $K(\cdot)$ given by expression \eqref{Kalpha}.

Since we want the above expression for any $X \sim S_{\alpha}(\sigma,\beta,\mu)$, we have the following result for $f_X(\cdot)$,
\begin{eqnarray}\label{fX1}
f_X (x)&=& \frac{1}{\sigma}f_{\alpha,\beta}\left(\frac{x-\mu}{\sigma} \right)\nonumber \\
&=& \frac{1}{\pi \sigma}\sum_{n=1}^\infty \frac{\Gamma(1+\frac{n}{\alpha})}{n!} \sin{\left(\pi n \frac{(\alpha + \beta'K(\alpha))}{2\alpha} \right)}\left(- \frac{x-\mu}{\sigma} \right)^{n-1} \nonumber \\
&=& \frac{1}{\pi \sigma}\sum_{n=1}^\infty \frac{\Gamma(1+\frac{n}{\alpha})}{n!} \sin{\left(\frac{\pi n}{2} \left(\frac{\alpha + \beta'(\alpha-2)}{\alpha}\right) \right)}\left(\frac{\mu-x}{\sigma} \right)^{n-1}.
\end{eqnarray}

We want the expression \eqref{fX1} for $\nu = n-1$ starting from zero. Hence,
\begin{equation}\label{fX2}
f_X (x)= \frac{1}{\pi \sigma}\sum_{\nu=0}^\infty \frac{\Gamma(1+\frac{\nu + 1}{\alpha})}{(\nu +1)!} \sin{\left(\frac{\pi (\nu +1)}{2\alpha}
\left(\alpha + \beta'\alpha-2\beta'\right) \right)}\left(\frac{\mu-x}{\sigma} \right)^{\nu}.
\end{equation}
Since $\Gamma(n)=(n-1)!$ and $\Gamma(z+1)= z\Gamma(z)$, expression \eqref{fX2} can be rewritten as
\begin{equation}\label{fX3}
f_X (x)= \frac{1}{\pi \sigma \alpha}\sum_{\nu=0}^\infty \frac{\Gamma(\frac{\nu + 1}{\alpha})}{\Gamma(\nu +1)} \sin{\left(\frac{\pi (\nu +1)}{2\alpha}
\left(\alpha + \beta'\alpha-2\beta'\right) \right)}\left(\frac{\mu-x}{\sigma} \right)^{\nu}, \ x \ge \mu.
\end{equation}

If $x<0$ we need to apply expression (2.7) in \cite{Schneider}, that is, $f_{\alpha,\beta}(-x)= f_{\alpha,-\beta}(x)$.
Therefore, expression \eqref{fX3} can be rewritten as
\begin{equation}\label{fX4}
f_X (x)= \frac{1}{\pi \sigma \alpha}\sum_{\nu=0}^\infty \frac{\Gamma(\frac{\nu + 1}{\alpha})}{\Gamma(\nu +1)} \sin{\left(\frac{\pi (\nu +1)}{2\alpha} \left(\alpha - \beta'\alpha+2\beta'\right) \right)}\!\left(\frac{x-\mu}{\sigma} \right)^{\nu}\! \!, \ x < \mu.
\end{equation}

Expressions \eqref{fX3} and \eqref{fX4} give us the central part of the density function $f_X (\cdot)$ for, respectively, $x \ge \mu$ and $x<\mu$.
Now we want the tails of $f_X (\cdot)$, when $\alpha >1$. We consider the expression (2.18) in \cite{Schneider},
\begin{equation*}
f_{\alpha, \beta} (x) \cong \frac{1}{\pi}\sum_{n=1}^\infty \frac{\Gamma(1+n\alpha)}{\Gamma{(n+1)}}(-1)^{n-1} \sin{(\pi n \alpha \gamma)}(x)^{-1-n\alpha}.
\end{equation*}

The equivalent expression for $f_X(\cdot)$, when $X \sim S_{\alpha}(\sigma,\beta,\mu)$, is given by
\[
\begin{array}{l}
f_X (x)\cong \\
\displaystyle \cong \frac{1}{\pi \sigma}\sum_{\nu=0}^\infty \frac{\Gamma(1+(\nu+1)\alpha)}{\Gamma{(\nu+2)}}(-1)^\nu \sin{\left(\pi (\nu+1) \alpha \left(\frac{\alpha + \beta'(\alpha-2)}{2\alpha}\right) \right)}\left(\frac{x-\mu}{\sigma} \right)^{-1-(\nu+1)\alpha} \nonumber \\
\displaystyle =\frac{1}{\pi \sigma}\sum_{\nu=0}^\infty \frac{(\nu+1)\alpha\Gamma((\nu+1)\alpha)}{(\nu+1)\Gamma{(\nu+1)}}(-1)^\nu
\sin{\!\left(\frac{\pi (\nu+1) \alpha}{2\alpha} \left(\alpha + \alpha\beta'-2\beta'\right) \right)} \! \! \left(\frac{x-\mu}{\sigma}
\right)^{-1-(\nu+1)\alpha} \nonumber \\
\displaystyle =\frac{\alpha}{\pi \sigma}\sum_{\nu=0}^\infty \frac{\Gamma((\nu+1)\alpha)}{\Gamma{(\nu+1)}}(-1)^\nu \left(\frac{\sigma}{x-\mu} \right) \sin{\!\left(\frac{\pi (\nu+1) }{2} \left(\alpha + \alpha\beta'-2\beta'\right) \right)} \! \! \left(\frac{\sigma}{x-\mu}\right)^{(\nu+1)\alpha}.
\nonumber
\end{array}
\]
Hence, for $x > \mu$,
\begin{eqnarray}\label{fX5}
f_X (x)\cong \frac{\alpha}{\pi}\sum_{\nu=1}^\infty \frac{\Gamma{((\nu+1)\alpha)}(-1)^\nu}{\Gamma{(\nu+1)}(x-\mu)}\sin{\!\left(\frac{\pi (\nu+1) }{2} \left(\alpha + \alpha\beta'-2\beta'\right) \right)} \! \! \left(\frac{\sigma}{x-\mu}\right)^{(\nu+1)\alpha}.
\end{eqnarray}

Again we apply expression (2.18) in \cite{Schneider} to obtain expression \eqref{fX5} above, when $x < \mu$
\begin{eqnarray}\label{fX6}
f_X (x)\cong \frac{\alpha}{\pi}\sum_{\nu=0}^\infty \frac{\Gamma{((\nu+1)\alpha)}(-1)^\nu}{\Gamma{(\nu+1)}(\mu-x)}\sin{\!\left(\frac{\pi (\nu+1) }{2} \left(\alpha - \alpha\beta' + 2\beta'\right) \right)} \! \! \left(\frac{\sigma}{\mu-x}\right)^{(\nu+1)\alpha}.
\end{eqnarray}

\vspace{0.3cm}

\noindent\textbf{Case 2: $\alpha <1$}

Now we want the formulas for $\alpha < 1$. Page 500 of \cite{Schneider} explains that the series expansion
of $f_{\alpha,\beta}(\cdot)$ is given by the right-hand side of (2.18) and its asymptotic behavior is given by the right-hand side of (2.17).

The series expansion of $f_X(\cdot)$, in the case $\alpha < 1$, where $K(\alpha)=\alpha$ is given by
\[
\begin{array}{l}
f_X (x)= \\
\displaystyle = \frac{1}{\pi \sigma}\sum_{n=1}^\infty \frac{\Gamma(1+n\alpha)}{n!}(-1)^{n-1} \sin{\left(\pi n\alpha\gamma \right)}
\left(\frac{x-\mu}{\sigma} \right)^{-1-n\alpha} \nonumber \\
\displaystyle = \frac{1}{\pi \sigma}\sum_{\nu=0}^\infty \frac{\Gamma(1+(\nu+1)\alpha)}{(\nu+1)!}(-1)^\nu \sin{\left(\pi (\nu+1) \alpha
\left(\frac{\alpha + \beta'\alpha}{2\alpha}\right) \right)}\left(\frac{x-\mu}{\sigma} \right)^{-1-(\nu+1)\alpha} \nonumber \\
\displaystyle =\frac{\alpha}{\pi}\sum_{\nu=0}^\infty \frac{\Gamma((\nu+1)\alpha)}{\Gamma{(\nu+1)}(x-\mu)}(-1)^\nu
\sin{\!\left(\frac{\pi (\nu+1) }{2} (\alpha + \alpha\beta') \right)} \! \! \left(\frac{\sigma}{x-\mu}\right)^{(\nu+1)\alpha}. \nonumber
\end{array}
\]
Hence, for $x > \mu$,
\begin{eqnarray}\label{fX7}
f_X (x)= \frac{\alpha}{\pi}\sum_{\nu=0}^\infty \frac{\Gamma((\nu+1)\alpha)}{\Gamma{(\nu+1)}(x-\mu)}(-1)^\nu \sin{\!\left(\frac{\pi (\nu+1) }{2}
(\alpha + \alpha\beta') \right)} \! \! \left(\frac{\sigma}{x-\mu}\right)^{(\nu+1)\alpha}.
\end{eqnarray}
The equivalent expression when $x<\mu$, is given by

\begin{eqnarray}\label{fX8}
f_X (x)= \frac{\alpha}{\pi}\sum_{\nu=0}^\infty \frac{\Gamma((\nu+1)\alpha)}{\Gamma{(\nu+1)}(\mu-x)}(-1)^\nu
\sin{\!\left(\frac{\pi (\nu+1) }{2} (\alpha - \alpha\beta') \right)} \! \! \left(\frac{\sigma}{\mu-x}\right)^{(\nu+1)\alpha}.
\end{eqnarray}

Now we shall give the formulas for the tails of the distribution when $\alpha <1$.
We consider the right-hand side of (2.17) in \cite{Schneider}.
The asymptotic behavior for $f_X(\cdot)$ is given by
\[
\begin{array}{rl}
f_X (x)\cong & \displaystyle \frac{1}{\pi \sigma\alpha}\sum_{\nu=0}^\infty \frac{\Gamma(\frac{\nu+1}{\alpha})}{\Gamma{(\nu+1)}}
\sin{\left(\frac{\pi (\nu+1)}{ 2\alpha} (\alpha + \alpha\beta')\right)}\left(\frac{x-\mu}{\sigma} \right)^\nu \nonumber \\
=& \displaystyle \frac{1}{\pi \sigma\alpha}\sum_{\nu=0}^\infty \frac{\Gamma(\frac{\nu+1}{\alpha})}{\Gamma{(\nu+1)}}
\sin{\left(\frac{\pi (\nu+1)}{ 2} (1 + \beta')\right)}\left(\frac{x-\mu}{\sigma} \right)^\nu. \nonumber
\end{array}
\]
Hence, the expression for the tails when $\alpha < 1$ is given by
\begin{eqnarray}\label{fX9}
f_X (x)\cong \frac{1}{\pi \sigma\alpha}\sum_{\nu=0}^\infty \frac{\Gamma(\frac{\nu+1}{\alpha})}{\Gamma{(\nu+1)}} \sin{\left(\frac{\pi (\nu+1)}{ 2} (1 + \beta')\right)}\left(\frac{x-\mu}{\sigma} \right)^\nu, \ x \ge \mu,
\end{eqnarray}
and its equivalent expression when $x<\mu$ is given by
\begin{eqnarray}\label{fX10}
f_X (x)\cong \frac{1}{\pi \sigma\alpha}\sum_{\nu=0}^\infty \frac{\Gamma(\frac{\nu+1}{\alpha})}{\Gamma{(\nu+1)}}
\sin{\left(\frac{\pi (\nu+1)}{2} (1 - \beta')\right)}\left(\frac{\mu - x}{\sigma} \right)^\nu, \ x < \mu.
\end{eqnarray}

We calculate the likelihood function based on the generated samples by considering expressions \eqref{fX3} and \eqref{fX4},
or \eqref{fX5} and \eqref{fX6}, for their central or tail components depending on the sign of $x- \mu$, when $\alpha>1$.
When $\alpha< 1$, we calculate the likelihood function based on the generated samples by considering expressions \eqref{fX7} and \eqref{fX8},
or \eqref{fX9} and \eqref{fX10}, for their central or tail components depending on the sign of $x- \mu$.
For the Cosine process, the noise
$\varepsilon_{k,h}$ is considered to be an independent identically distributed sequence of random variables $S_{\alpha}(\sigma_\varepsilon,0,0)$, with $\sigma_\varepsilon^{\alpha}=2\int_0^h|\cos(as)|^{\alpha}ds$.
Let $\mathbf{V}=\{V_{kh}\}_{k=0}^{n-1}$ be a sample from the Cosine process and $\bm{\eta}=(\alpha,\sigma_\varepsilon,a)'$ be the parameter
vector to be estimated. Let $V_0$ and $V_1$ be two random variables $S_{\alpha}(\sigma_0,0,0)$ distributed with $\sigma_0=1.0$.
Then, the likelihood function is given by
\begin{align*}
\mathcal{L}(\mathbf{V}|\bm{\eta})=&\prod_{k=0}^{n-1} f(V_{(k+1)h}-2\cos(ah)V_{kh}+V_{(k-1)h};\sigma_\varepsilon,0,0),
\end{align*}
\noindent where $f(\cdot)$ is given by the expressions \eqref{fX3}-\eqref{fX10}, depending on
the values of $\alpha$ and on the sign of $x - \mu$.

By considering the independence of the parameters, we set the following non-informative \emph{priori} distributions:

\vspace{-0.3cm}
\begin{itemize}
\item $\alpha \sim {\cal U}([0,2])$, that is, $\alpha$ follows a Uniform distribution on the interval $[0,2]$;

\item $a \sim {\cal U}([0,3])$, that is, $a$ follows a Uniform distribution on the interval $[0,3]$;

\item $\sigma_{\varepsilon} \sim \Gamma(1,2)$, that is, $\sigma_{\varepsilon}$ follows a Gamma distribution with parameters 1 and 2.
\end{itemize}
\vspace{-0.3cm}

Therefore, the posterior distribution is given by
\begin{align*}
\Pi(\bm{\eta}|\mathbf{V})&\propto \mathcal{L}(\mathbf{V}|\bm{\eta})\pi(\bm{\eta}),
\end{align*}
\noindent where $\pi(\cdot)$ is the priori distribution of the parameter vector $\bm{\eta}=(\alpha, \sigma_{\varepsilon}, a)'$.

We consider the following scenarios for the Bayesian estimation procedure:
\[
\bullet \ \ a\in\{1,2\};  \ \ \bullet \ \ \alpha\in\{1.1,1.5,2\}; \ \  \bullet \ \ n=2000; \  \  \bullet \ \ h=1.
\]

We generate $30,000$ Gibbs samples for each one of the three parameters. We use a burn-in-sample of size $10,000$ and we take every $10$-th sample which  gives a final sample of size $2,000$ to be used for finding the posterior summaries of interest. Table \ref{tabelaCOSbayesiano.estavel} presents the Bayesian method results. For each parameter in the vector $\bm{\eta}=(\alpha, \sigma_{\varepsilon}, a)'$, Table \ref{tabelaCOSbayesiano.estavel} reports the mean, the bias, the standard error (s.e.) and the credibility interval (CI) for its estimator. We notice that these results present low bias values for all the estimators. However, the true $\alpha$ parameter value is not inside its credibility interval when $\alpha \in \{1.1, 2\}$.
The $a$ parameter is estimated in a very accurate way. This is true for all scenarios considered here.

Similarly to the OU process, for the Cosine process we also have problems related to the initial values of the chains.To mitigate these problems we considered an approximated version of the estimated means obtained from the maximum likelihood procedure to initiate the chains.

\begin{figure}[h!!!]
\begin{footnotesize}
\centering

\mbox{
   \subfigure[Traceplot for $\alpha$]{\includegraphics[width=0.25\textwidth]{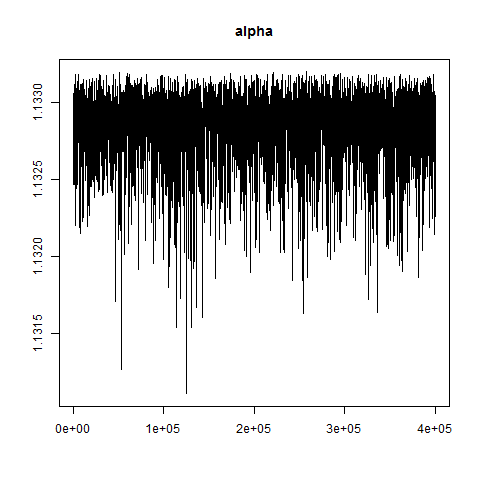}}
   \subfigure[Density function $\alpha$]{\includegraphics[width=0.25\textwidth]{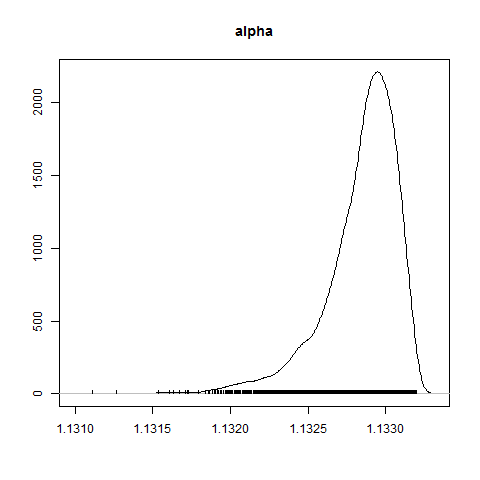}}
   \subfigure[acf $\alpha$]{\includegraphics[width=0.25\textwidth]{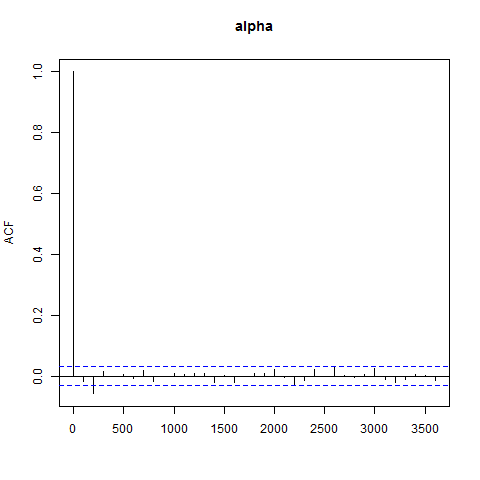}}
  }
\mbox{
   \subfigure[Traceplot for $\sigma_{\varepsilon}$]{\includegraphics[width=0.25\textwidth]{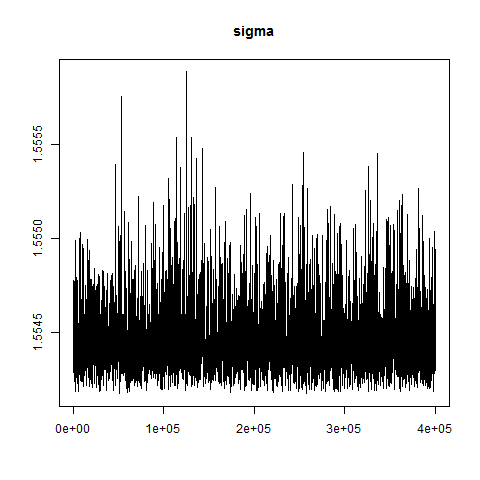}}
   \subfigure[Density function $\sigma_{\varepsilon}$]{\includegraphics[width=0.25\textwidth]{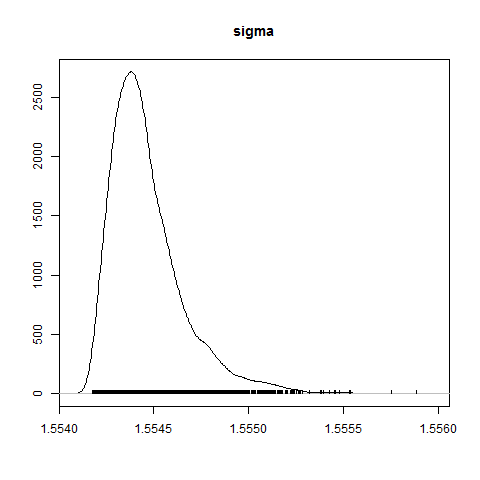}}
   \subfigure[acf $\sigma_{\varepsilon}$]{\includegraphics[width=0.25\textwidth]{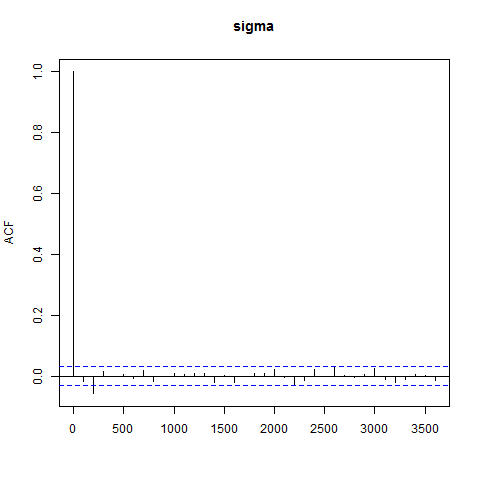}}
  }
  \mbox{
   \subfigure[Traceplot for $a$]{\includegraphics[width=0.25\textwidth]{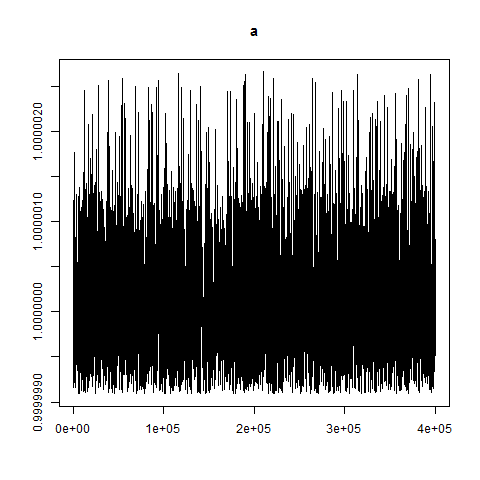}}
   \subfigure[Density function $a$]{\includegraphics[width=0.25\textwidth]{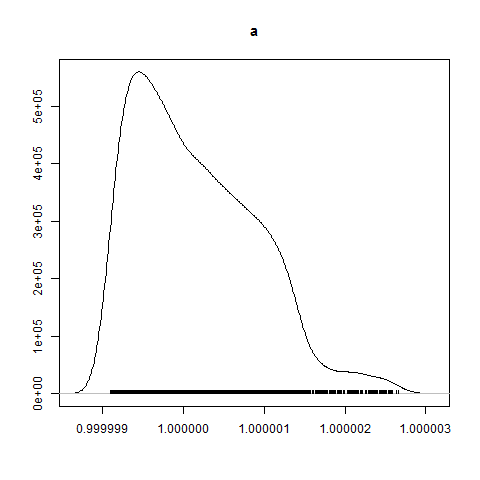}}
   \subfigure[acf $a$]{\includegraphics[width=0.25\textwidth]{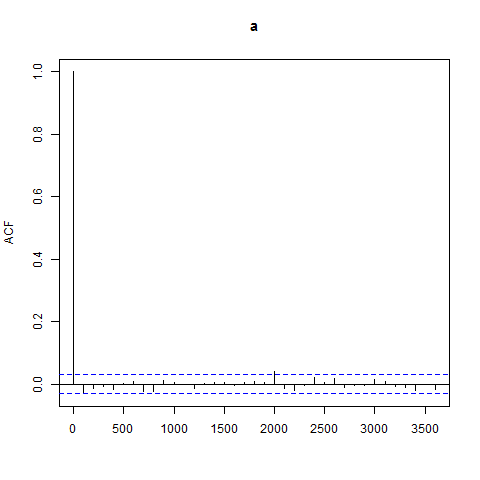}}
  }
  \caption{Traceplot for the generated chains, density and autocorrelation (acf) functions, when $a=1$, $\alpha=1.1$, $n=2000$ and $h=1$.}\label{figCOSbayesiano.estavel}
  \end{footnotesize}
\end{figure}

Figure \ref{figCOSbayesiano.estavel} presents the traceplot for the generated chains, the density and the autocorrelation functions
for the case when $a=1$, $\alpha=1.1$, $n=2000$ and $h=1$. All other cases had similar performance and we omitted to report them here.

\begin{table}[!h]
\caption{Bayesian estimation results for the Cosine process, when $a\in\{1,2\}$, $\alpha\in\{1.1,1.5,2\}$, $n=2000$ and $h=1$.}\label{tabelaCOSbayesiano.estavel}
\vspace{-0.5cm}
\begin{center}
{\footnotesize
\begin{tabular}{|c|ccc|}
\hline
\hline
Statistic   & $\alpha=1.1$  & $\sigma_\varepsilon=1.5824$ & $a=1$     \\
\hline
\hline

      $\mbox{mean}$ &     1.1787 &     1.6296 &     0.9999 \\

      $\mbox{bias}$ &    0.0070 &   0.0469  &    0.0001\\

      $\mbox{s.e.}$ &    0.0002 &     0.0010   &     0.0000 \\

      $\mbox{CI}$ &    [1.1713,1.1964] &  [1.5400,1.7232]  & [0.9997,1.0000] \\
\hline
\hline
Statistic    & $\alpha=1.1$ & $\sigma_\varepsilon=1.0452$ & $a=2$       \\
\hline
      $\mbox{mean}$ &     1.1782 &      1.0792  &    2.0000 \\

      $\mbox{bias}$ &    0.0069 &   0.0623  &    0.0001\\

      $\mbox{s.e.}$ &    0.0002 &     0.0014  &     0.0000  \\

      $\mbox{CI}$ &    [1.1713,1.1965] &  [1.0189,1.1387]  & [1.9998,2.0003] \\
\hline
\hline
Statistic    & $\alpha=1.5$ & $\sigma_\varepsilon=1.3450$ & $a=1$    \\
\hline
      $\mbox{mean}$ &     1.5072 &      1.3747 &    1.0001 \\

      $\mbox{bias}$ &    0.0335 &   0.0340 &    0.0003\\

      $\mbox{s.e.}$ &    0.0007 &     0.0008  &     0.0000  \\

      $\mbox{CI}$ &    [1.4401,1.5701] &  [1.3088,1.4415]  & [0.9996,1.0006]   \\
\hline
\hline
Statistic   & $\alpha=1.5$ & $\sigma_\varepsilon=0.9467$ & $a=2$    \\
\hline
      $\mbox{mean}$ &    1.5082 &      0.9679 &   2.0002 \\

      $\mbox{bias}$ &    0.0325 &   0.0229 &    0.0003\\

      $\mbox{s.e.}$ &   0.0007 &     0.0005  &    0.0000 \\

      $\mbox{CI}$ &    [1.4420,1.5699] &  [0.9249,1.0143]  & [1.9995,2.0009] \\
\hline
\hline
Statistic & $\alpha=2$ & $\sigma_\varepsilon=1.2061$ & $a=1$    \\
\hline
      $\mbox{mean}$ &     1.9927 &      1.2115 &    1.0004 \\

      $\mbox{bias}$ &    0.0069 &   0.0198 &    0.0011\\

      $\mbox{s.e.}$ &    0.0002 &     0.0004  &     0.0000  \\

      $\mbox{CI}$ &    [1.9737,1.9998] &  [1.1739,1.2505]  & [0.9983,1.0025]  \\
\hline
\hline
Statistic   & $\alpha=2$ & $\sigma_\varepsilon=0.9004$ & $a=2$    \\
\hline
      $\mbox{mean}$ &     1.9930 &      0.9043 &    2.0010 \\

      $\mbox{bias}$ &    0.0066 &   0.0148 &    0.0006\\

      $\mbox{s.e.}$ &    0.0001 &     0.0003   &     0.0000  \\

      $\mbox{CI}$ &    [1.9755,1.9999] &  [0.8761,0.9334]  & [1.9999,2.0021] \\
\hline
\hline
\end{tabular}}
\end{center}
\end{table}

\begin{rmk}
Although the computational time for the Bayesian estimation MCMC procedure, based on the MCMC algorithm, is much higher than for the maximum likelihood method, we would like to mention some advantages:

\begin{itemize}
\item the possibility of selecting the prior distribution based on a piece of known information from the literature or provided by a specialist;

\item bias and standard error values for the parameters tend to be smaller when the Bayesian approach is considered;

\item a sample from the posterior distribution of each parameter is obtained, rather than a single point, which makes it possible, if one wishes, to construct confidence intervals and/or perform a hypothesis test without making further assumptions on the distribution of the estimates.

\end{itemize}

However, in our simulation study, we found a disadvantage for the used Bayesian estimation procedure: the credibility intervals do not always
contain the true parameter value for the stability index $\alpha$, when it is in the set $\{1.1, 2\}$. However, for the true parameters $a$
and $\sigma_{\varepsilon}$, they are always inside the credibility interval for all cases considered here.

\end{rmk}

\subsection{Goodness-of-fit Tests}

We present some goodness-of-fit tests that shall be used in the
Applications Section. Goodness-of-fit tests are important devices to decide
if an observed data set follows a specific distribution or not.
It might appear that a certain real data set is better fitted by a Gaussian distribution ($\alpha=2$), when, in fact, an $\alpha$-stable
non-Gaussian distribution ($0 <\alpha< 2$) does better.
We know that $\alpha$-stable distributions may represent a good fit for modeling
the financial asset returns and that Normality is not a good choice (see \cite{Borak}). In Application 1 we will see this fact and the \emph{Cosine Process} driven by an $\alpha$-stable noise is the best model to fit the real data set. However, in
Application 2, we will show that the \emph{Cosine Process} driven by a Gaussian noise best fitted this real data set.

We shall consider the classical Kolmogorov-Smirnov (KS) and Anderson-Darling (AD)
tests (see \cite{DAgostino}). Besides these two goodness-of-fit tests we shall also consider a modified version of the Kolmogorov-Smirnov (MKS) (see \cite{Beaulieu}) and
the McCulloch (MC) test, based on the estimators of \cite{McCulloch}.
The three former tests are known as \emph{empirical distribution functions} (EDF) tests
since they all rely on some distance between the empirical and the theoretical distribution functions. Moreover, the AD and MKS tests serve as an alternative to the classical KS test, since they give more weights to the tails than the original one (see \cite{Thode} and \cite{Beaulieu}).

Let us consider the ordered sample $\bm{X} = (X_{(1)}, \cdots, X_{(n)})$. The classical approach of goodness-of-fit tests consists on testing the null
hypothesis $H_0$ against the alternative $H_1$, given by
\begin{equation}\label{hypotheses1}
H_0: F_n(x) = S(x)  \ \ \mbox{vs} \ \ H_1:  F_n(x) \ne S(x),
\end{equation}
\noindent where $F_n(\cdot)$ is the empirical distribution function while
$S(\cdot)$ denotes the hypothetical distribution function.
When testing for Gaussian distributions, the standard procedure is to consider
$S(\cdot) = \Phi(\cdot)$, where $\Phi(\cdot)$ denotes the cumulative distribution
function of the standard Gaussian distribution ${\cal N}(0,1)$. If
\begin{align*}
 \bm{\bar{X}} = \frac{1}{n} \sum_{j=1}^n X_j \quad \mbox{and} \quad \bm{s_{X}} = \sqrt{\frac{\sum_{j=1}^n(X_j -\bm{\bar{X}})^2}{n-1}}
\end{align*}
are, respectively, the mean and standard deviation of the ordered observations $\bm{X}$, let $\bm{Y} = (Y_{(1)}, \cdots, Y_{(n)})$ be the standardized observations, given by
\begin{equation}\label{standard}
Y_{(j)} = \frac{X_{(j)} - \bm{\bar{X}}}{\bm{s_{X}}}, \quad \mbox{for } j = 1, \cdots, n.
\end{equation}

As \cite{Beaulieu} noticed, there is no need to consider the parameters $\sigma$ and $\mu$ when testing for $\alpha$-stable distributions. Indeed, if $X \sim S_\alpha(\sigma,\beta, \mu)$ then $Y = \frac{X - \mu}{\sigma}$ is an $\alpha$-stable
$Y \sim S_\alpha(1,\beta, 0)$ random variable.

In Section 6, we shall use the KS, AD, MKS and MC tests based on the hypotheses given by equation \eqref{hypotheses1}, for testing $\alpha$-stable distributions when real data observations are of interest. Since we are considering only symmetric $\alpha$-stable random
variables, we have $\beta=0$ and the hypothetical value
for $\hat{\alpha}$ is obtained by the maximum likelihood estimation procedure described in Section 5.1. Let $F_j^{n} = F_n(Y_{(j)})$ and $\hat{F}_j = F(Y_{(j)}; \hat{\alpha}),$ where $F_n(\cdot)$ is the empirical distribution function of the standardized data $\bm{Y}$,\ given in equation \eqref{standard}, and $F(\cdot; \alpha)$ is
the distribution function of a standard symmetric $\alpha$-stable distribution ($\sigma = 1$, $\beta=0= \mu$) with stability index $\alpha$. The goodness-of-fit statistical tests considered here are the following:
\begin{itemize}
 \item the KS statistic
\begin{align}\label{ks}
 D_n = \max_{j\,\in\,\{1,\cdots,n\}}\left[\hat{F}_j - F_{j-1}^{\,n},\ F_j^{\,n} - \hat{F}_j\right];
\end{align}
\item the AD statistic
\begin{align}\label{ad}
 A^2 =  -n -\sum_{j=1}^n \frac{(2j-1)}{n} \left[\log(\hat{F}_j) + \log(1 - \hat{F}_{n-j+1})\right];
\end{align}
\item the MKS statistic
\begin{align}\label{mks}
 MKS = \sqrt{n} \max_{j\,\in\,\{1,\cdots,n\}} \left[\frac{\left| F_j^{\,n} - \hat{F}_j \right|}{\hat{F}_j \left(1 - \hat{F}_j\right) + 1/n} \right];
\end{align}

\item the MC statistic
\begin{align}\label{MC}
\hat{\phi}_i (\alpha_0)= |\phi_i - \hat{\phi}_i (\alpha_0)|, \ \ \mbox{for} \ \ i=1,2,
\end{align}
where $\hat{\phi}_1$ and $\hat{\phi}_2$ are respectively the simulated values for $\phi_1$ and $\phi_2$ given by

\begin{align}
\phi_1= \frac{Y[95] - Y[5]}{Y[75] - Y[25]} \ \ \mbox{and} \ \
\phi_2=\frac{Y[95] + Y[5]- 2Y[50]}{Y[95] - Y[5]},
\end{align}
and $Y[z]$ denotes the $z$-th quantile of the standardized data given in
\eqref{standard}.

\end{itemize}

\section{Applications}
\renewcommand{\theequation}{\thesection.\arabic{equation}}
\setcounter{equation}{0}

In this section, we consider two real data sets. The first one deals with
the Apple company stock market price data while the second
considers the cardiovascular mortality in Los Angeles County data.

\subsection{Apple Company Stock Market Price Data}
In this subsection we analyze the behavior of the Apple company stock market price data. The period of this time series is from January 04, 2010, to November 26, 2014. The trading hours are from 9:31 am to 4:00 pm (Chicago time) and the tick-by-tick frequency of the original time series is one-minute, which gives a total of $390$ observations per day in a total of $1,235$ days. The overall total is $481,650$ observations.

For the analysis we consider four aggregated log-returns: the one-minute, the five-minute, the ten-minute and the fifteen-minute aggregated to obtain the four time series denoted, respectively, by $R_t^{(1)}$, $R_t^{(5)}$, $R_t^{(10)}$ and $R_t^{(15)}$. The one-minute log-returns are defined by
\begin{equation*}\label{logreturn}
R_t^{(1)}=100\times\ln\left(\frac{X_{t+1}}{X_{t}}\right),
\end{equation*}
\noindent where $t\in\{1,\cdots,481,649\}$ and $X_t$ is the price of the stocket market of Apple company at time $t$.

Figure \ref{originalandlogreturns} below shows both the original time series and its one-minute log-returns.

Based on the one-minute log-returns we aggregated them to obtain the five-minute, ten-minute, and fifteen-minute log-returns. Thus, we have the following four time series
\begin{itemize}
\item $R_t^{(1)}$, for $t\in\{1,\cdots,481,649\}$, with discretization step size given by $h=(390)^{-1}= 2.5641 \times 10^{-3}$;
\item $R_t^{(5)}= R_{5t-4}^{(1)}+R_{5t-3}^{(1)}+\cdots+R_{5t}^{(1)}$, for $t\in\{1,\cdots, 96,329\}$, with discretization
step
\newline
size given by $h=5(390)^{-1}= 0.0128$;
\item $R_t^{(10)}= R_{10t-9}^{(1)}+R_{10t-8}^{(1)}+\cdots+R_{10t}^{(1)}$, for $t\in\{1,\cdots, 48,164\}$, with discretization
step
\newline
size given by $h=10(390)^{-1}= 0.0256$;
\item $R_t^{(15)}= R_{15t-14}^{(1)}+R_{15t-13}^{(1)}+\cdots+R_{15t}^{(1)}$, for $t\in\{1,\cdots,32,109\}$, with discretization
step size given by $h=15(390)^{-1}= 0.0385$.
\end{itemize}

\begin{figure}[H]
\centering
\mbox{
\subfigure[]{\includegraphics[width=0.48\textwidth]{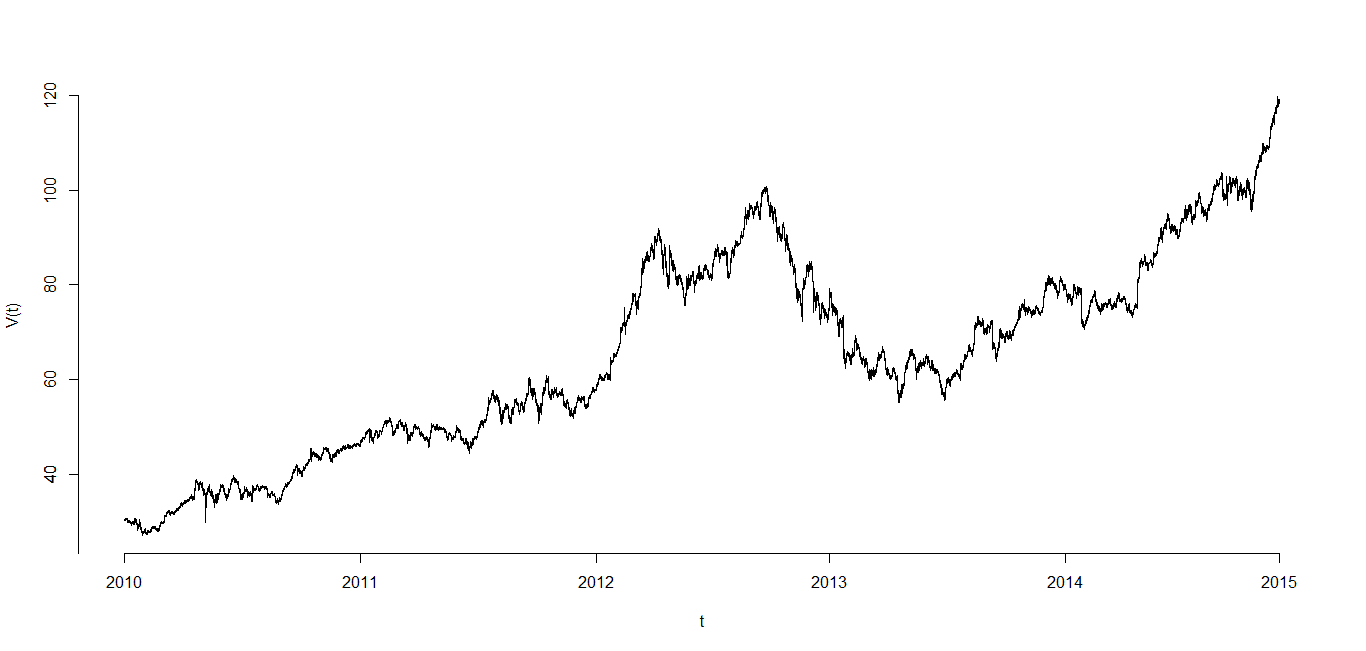}}
\subfigure[]{\includegraphics[width=0.48\textwidth]{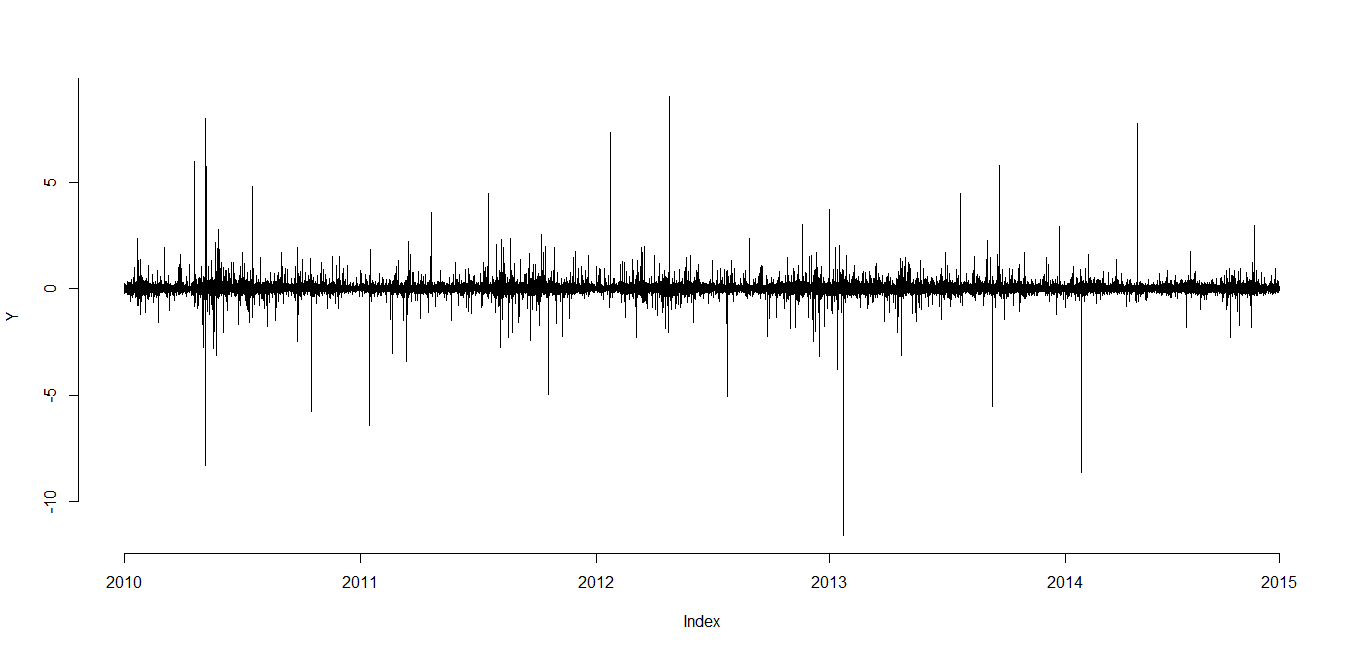}}
}
\caption{Stock market price of Apple company, from January 04, 2010 to November 26, 2014:
(a) original time series; (b) one-minute log-return time series.}
\label{originalandlogreturns}
\end{figure}

\vspace{-0.7cm}

\begin{figure}[h!!!]
\centering
\mbox{
\subfigure[]{\includegraphics[width=0.48\textwidth]{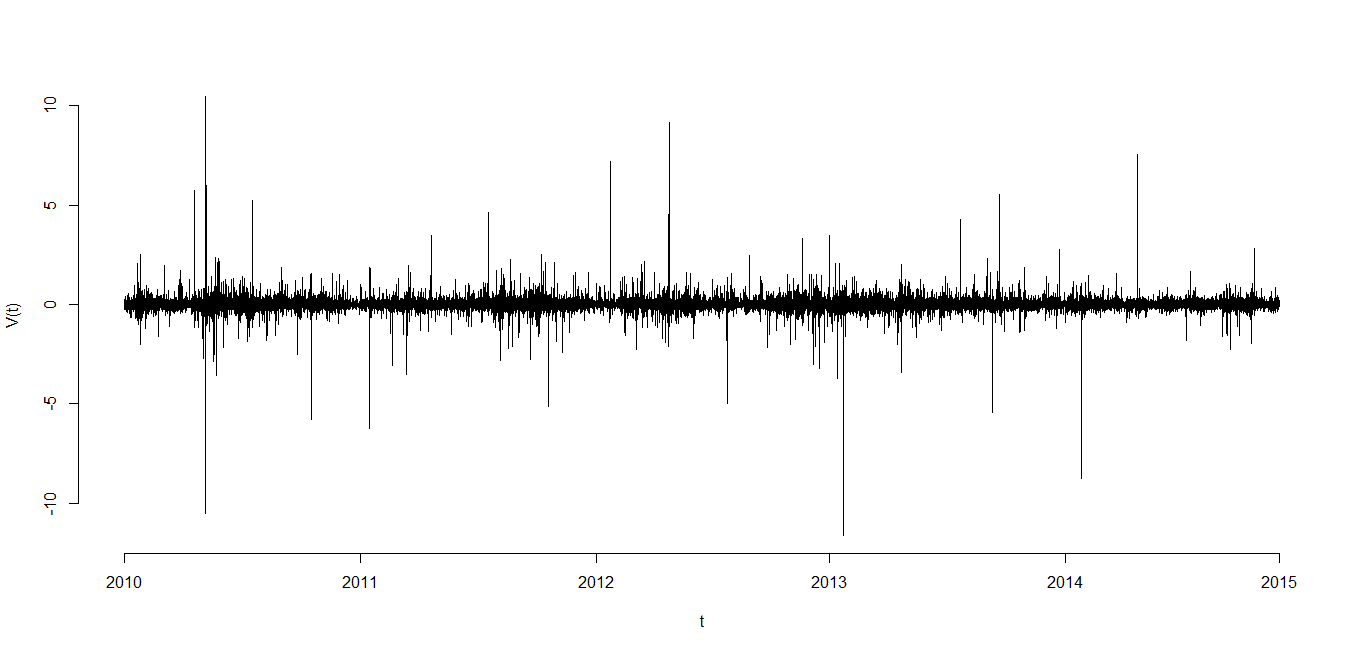}}
\subfigure[]{\includegraphics[width=0.48\textwidth]{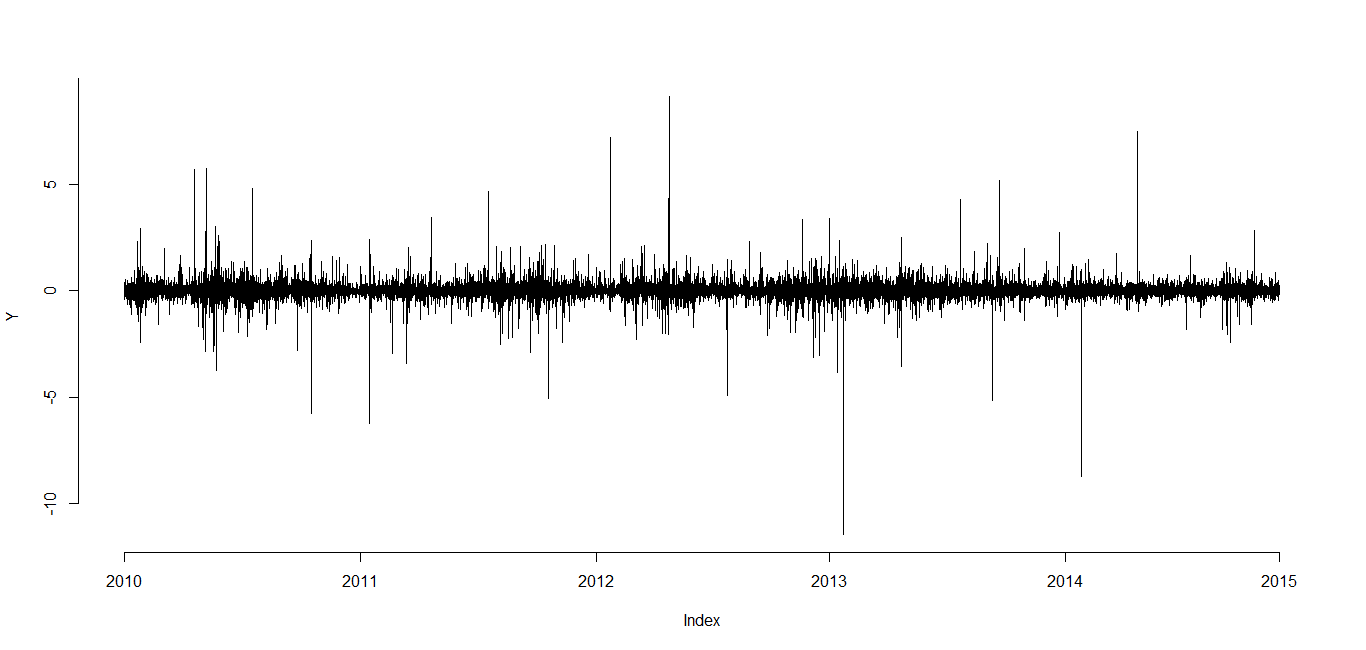}}
}
\mbox{
\subfigure[]{\includegraphics[width=0.48\textwidth]{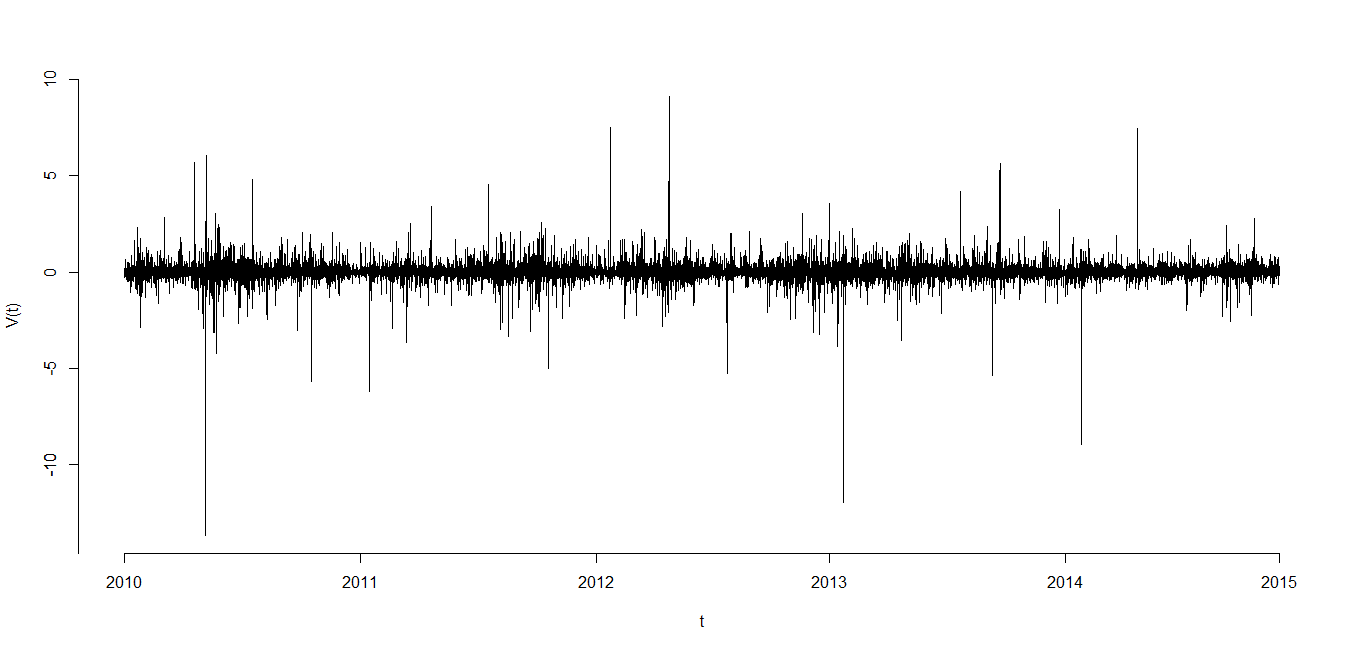}}
}
\caption{Stock market price of Apple company: (a) five-minute log-returns; (b) ten-minute log retorns; (c) fifteen-minute log-returns.}
\label{fiveandfifteenlogreturns}
\end{figure}

\begin{figure}[H]
\centering
\mbox{
\subfigure[]{\includegraphics[width=0.24\textwidth]{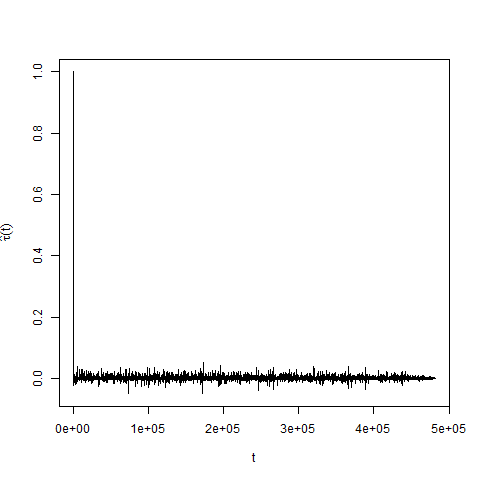}}
\subfigure[]{\includegraphics[width=0.24\textwidth]{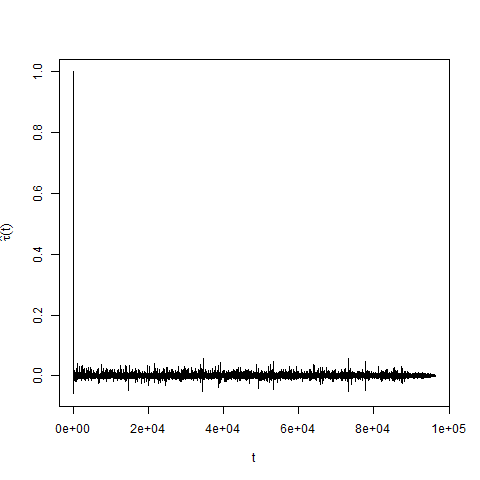}}}
\subfigure[]{\includegraphics[width=0.24\textwidth]{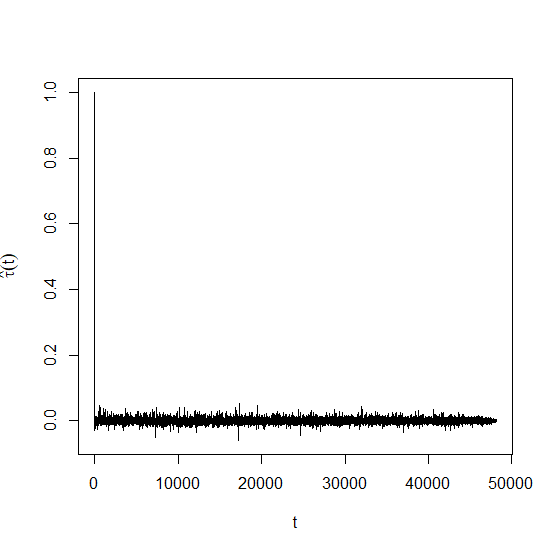}}
 \subfigure[]{\includegraphics[width=0.24\textwidth]{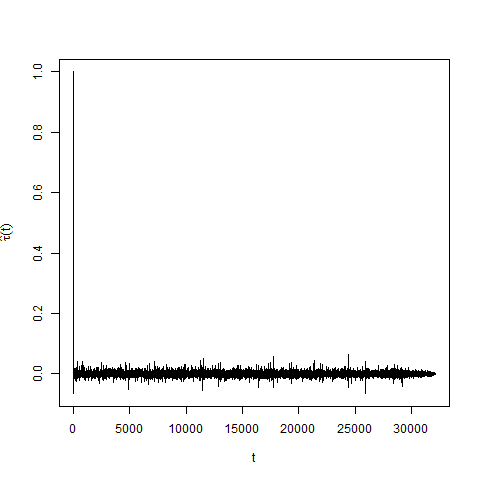}}

\vspace{-0.3cm}
\caption{Normalized empirical codifference function, when $s=0.01$, for the: (a) one-minute log-return; (b) five-minute log-returns; (c) ten-minute log-returns; (d) fifteen-minute log-returns.}
\label{codifApple}
\end{figure}

Figure \ref{fiveandfifteenlogreturns} shows the three time series with, respectively, five-minute, ten-minute and fifteen-minute aggregated log-returns
of the stock market price for the Apple company. Figure \ref{codifApple}  shows the normalized empirical codifference for the one-minute,
five-minute, ten-minute and fifteen-minute log-return time series when $s= 0.01$.

To model these four log-return time series denoted, respectively, by $R_t^{(1)}$, $R_t^{(5)}$, $R_t^{(10)}$ and $R_t^{(15)}$,
we consider the Cosine process, given in expression \eqref{eqSOL}. In this case, the memory function is given by $\rho(t) = \cos(at)$, for all $t\geq 0$ where $a$ is a positive real constant. We consider the maximum likelihood estimation procedure, presented in Section 5.1, to estimate the
vector parameter $\bm{\eta}=(\alpha, \sigma_{\varepsilon}, a)$. Table \ref{AppleMLE} presents the estimation results
considering the four log-return time series, $R_t^{(j)}$, based
on the four respective discretization steps with sizes $h = j \, (390)^{-1}$, for $j \in \{1,5,10,15\}$. We know that log-return
time series are well modeled by heavy tailed processes (see, for instance, \cite{Lopes2013}, \cite{Lopes2014} and \cite{Prass}).

\begin{figure}[H]
\centering
\mbox{
\subfigure[$R_t^{(1)}$; $h=0.0026$]{\includegraphics[width=0.33\textwidth]{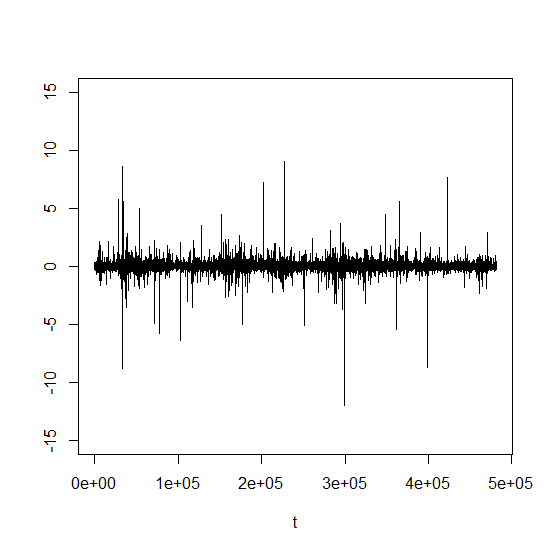}} \subfigure[$R_t^{(5)}$; $h=0.0128$]{\includegraphics[width=0.33\textwidth]{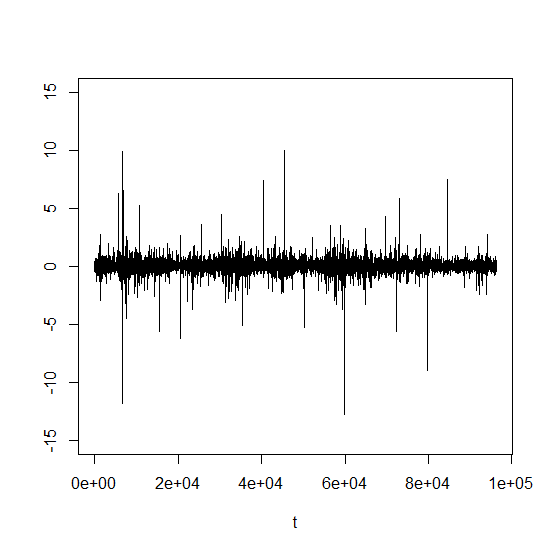}}}
\mbox{
\subfigure[$R_t^{(10)}$; $h=0.0256$]{\includegraphics[width=0.33\textwidth]{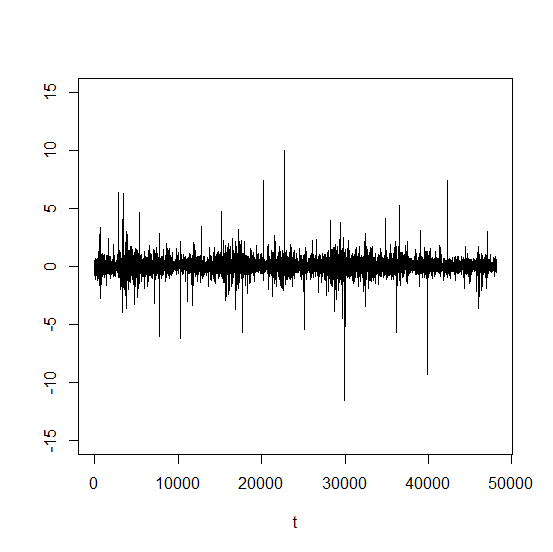}}
 \subfigure[$R_t^{(15)}$; $h=0.0385$]{\includegraphics[width=0.33\textwidth]{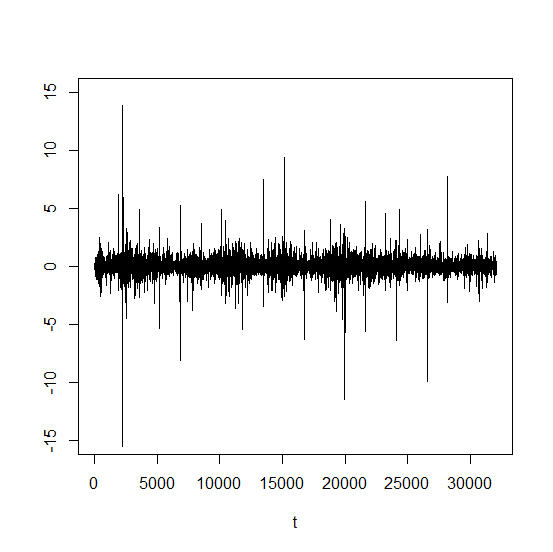}}
}
\vspace{-0.3cm}
\caption{Residual values for the four fitted log-return time series: (a) $R_t^{(1)}$; (b) $R_t^{(5)}$; (c) $R_t^{(10)}$; (d) $R_t^{(15)}$.}
\label{residuosApple}
\end{figure}

\begin{table}[h!!!]\caption{Maximum likelihood estimation results, based on the aggregated log-return time series $R_t^{(j)}$, with respective discretization step size $h =j\,(390)^{-1}$, for $j \in \{1,5,10,15\}$.}\label{AppleMLE}
\centering
\begin{tabular}{|c|c|c|}
\hline
\hline
   $\hat{\alpha}$ & $\hat{\sigma}_{\varepsilon}$ &  $\hat{a}$ \Tspace\Bspace\\
\hline
\hline
\multicolumn{ 3}{|c|}{$R_t^{(1)}$; $h=(390)^{-1}= 2.564 \times 10^{-3}$}  \Tspace\Bspace\\
\hline
1.5132   & 0.0472   &   632.5659       \Tspace\Bspace\\
\hline
\hline
\multicolumn{ 3}{|c|}{$R_t^{(5)}$; $h=5(390)^{-1}=0.0128$}  \Tspace\Bspace\\
\hline
1.4661   &    0.1010   & 124.6056   \Tspace\Bspace\\
\hline
\hline
\multicolumn{ 3}{|c|}{$R_t^{(10)}$; $h=10(390)^{-1}=0.0256$}  \Tspace\Bspace\\
\hline
1.4531   &    0.1468   & 62.9701   \Tspace\Bspace\\
\hline
\hline
\multicolumn{ 3}{|c|}{$R_t^{(15)}$; $h=15(390)^{-1}=0.0385$}  \Tspace\Bspace\\
\hline
1.4304  &  0.1959   &  36.6624    \Tspace\Bspace\\
\hline
\hline
\end{tabular}
\end{table}

\begin{table}[H]\caption{Statistics values for the Kolmogorov-Smirnov (KS),
Anderson-Darling (AD) and Modified Kolmogorov-Smirnov
(MKS) tests, applied to four aggregated log returns of the Apple company stock market price,
with their corresponding $p$-values (in parentheses).}\label{TestsApplic1}
\centering
\begin{tabular}{ccc}
\hline
\hline
\mbox{KS} & \mbox{AD} & \mbox{MKS} \\
\hline
\hline
\multicolumn{3}{c}{$R_t^{(1)}$: $H_0$: $\alpha=2$} \\
\hline
0.1928(0.001) & 48449.56(0.001)& 239122.90(0.001) \\
\hline
\multicolumn{3}{c}{$R_t^{(1)}$: $H_0$: $\hat{\alpha}=1.5132$} \\
\hline
0.1927(1) & 49709.4 (1) & 689.49(1)  \\
\hline
\hline
\multicolumn{3}{c}{$R_t^{(5)}$: $H_0$: $\alpha=2$} \\
\hline
0.1976(0.001) & 9969.62(0.001) & 24265.75(0.001)    \\
\hline
\multicolumn{3}{c}{$R_t^{(5)}$: $H_0$: $\hat{\alpha}=1.4661$} \\
\hline
0.1970(1) & 10243.03(1) & 306.73(1)    \\
\hline
\hline
\multicolumn{3}{c}{$R_t^{(10)}$: $H_0$: $\alpha=2$} \\
\hline
0.1840(0.001) & 4344.31(0.001) & 10202.57(0.001)     \\
\hline
\multicolumn{3}{c}{$R_t^{(10)}$: $H_0$: $\hat{\alpha}=1.4531$} \\
\hline
0.1836(1) & 4514.39(1) & 217.47(1)    \\
\hline
\hline
\multicolumn{3}{c}{$R_t^{(15)}$: $H_0$: $\alpha=2$} \\
\hline
0.1903(0.001) & 3084.83(0.001) & 5406.17(0.001)    \\
\hline
\multicolumn{3}{c}{$R_t^{(15)}$: $H_0$: $\hat{\alpha}=1.4304$} \\
\hline
0.1895(1)& 3214.16(1) & 175.30(1)    \\
\hline
\hline
\end{tabular}
\end{table}

Figure \ref{residuosApple} shows the residuals obtained from the fitted Cosine model, for all four considered log-return time series.
From the four panels of this figure, we can see that the Cosine process based on the original log-return time series $R_t^{(1)}$ captures
well the features of this data set.  However, the one-minute log return time series with
the smallest discretization step size also has the
smallest standard deviation estimate and the smallest variation of the residual values. Table \ref{AppleMLE} presents the maximum likelihood estimation for the parameters of
all four log return time series.

Table \ref{TestsApplic1} presents the EDF tests given in Section 5.3
for the four log return time series. The statistical tests and their
respective $p$-value (in parentheses) for testing hypotheses $H_0$ versus $H_1$, defined in expression \eqref{hypotheses1} are in this table. For the EDF tests, the empirical distribution function is compared with the Gaussian one. This table also presents the corresponding statistics values when one considers the second hypotheses test defined as
\begin{equation}\label{hypotheses2}
H_0: \hat \alpha= \alpha_0 \ \ \mbox{vs} \ \ H_1: \hat \alpha
\neq \alpha_0,
\end{equation}
\noindent where $\hat{\alpha}$ is the maximum likelihood estimate value, obtained from the procedure described in Section 5.1, for each log return time series.
From the results in Table \ref{TestsApplic1} one observes, from the two hypotheses tests,
the four log return time series are considered to be best fitted by an $\alpha$-stable
distribution with their respective $\hat{\alpha}$ given in Table \ref{AppleMLE}.

\subsection{Cardiovascular Mortality in Los Angeles County}

For this second data set, we consider the cardiovascular mortality in
Los Angeles County (see \cite{Shumway}). This is a weekly time series from 1970 to
1979, with $n=508$ observations. Since the period of
this time series is $52$ weeks over the 10-year period, its frequency is given by
$a=\frac{2\pi}{52}=0.1208$. Figure 6.5 shows the strong seasonal component in
this series, corresponding to winter-summer variations and the downward
trend in cardiovascular mortality.

In the analysis, we shall consider the Cosine process to adjust this data set. For
the estimation procedure we consider the original time series $X(\cdot)$ and $Y(\cdot)$ resulting from two differentiations applied to the original one:
\begin{equation}
Z(t)=X(t+1)-X(t)\ \ \mbox{ and} \ \ Y(t)=Z(t+52)-Z(t).
\end{equation}

\begin{figure}[H]
\centering
\mbox{
\subfigure[]{\includegraphics[width=0.5\textwidth]{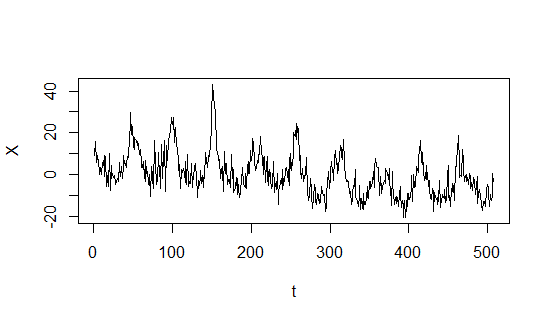}}
\subfigure[]{\includegraphics[width=0.5\textwidth]{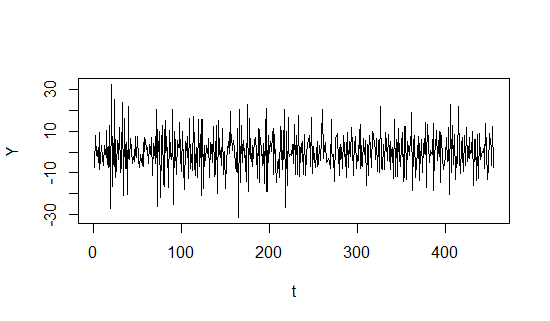}}}
\mbox{
\subfigure[]{\includegraphics[width=0.5\textwidth]{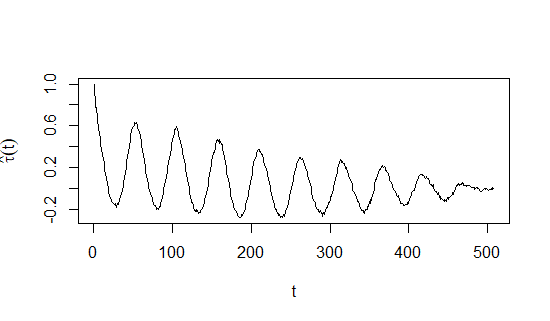}}
\subfigure[]{\includegraphics[width=0.5\textwidth]{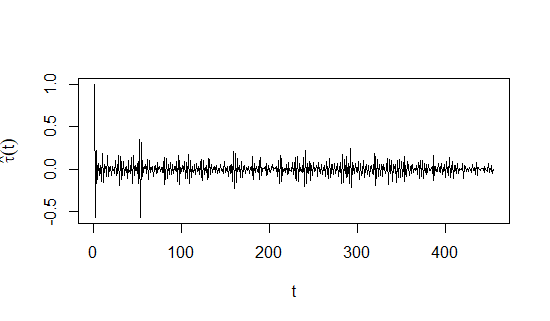}}}
\vspace{-0.3cm}
\caption{Cardiovascular mortality in Los Angeles County: (a) Original time series $X(\cdot)$; (b) Twice-differenced time series $Y(\cdot)$; (c) Normalized empirical codifference function for $X(\cdot)$; (d) Normalized empirical codifference function for $Y(\cdot)$.}
\label{figcardionormcodif}
\end{figure}

Figure \ref{figcardionormcodif}(a) shows the original time series $X(\cdot)$, centered
in its mean value while Figure \ref{figcardionormcodif}(b) shows the twice-differenced time series $Y(\cdot)$ without the strong seasonality and the downward trend. Figures \ref{figcardionormcodif}(c)-(d) present the normalized codifference for both time series
$X(\cdot)$ and $Y(\cdot)$. From Figure \ref{figcardionormcodif}(c) we see a
strong seasonal component as well as the downward trend in the original time series.

Table \ref{MortalityMLE} presents the maximum likelihood estimation for the parameters, for both times series $X(\cdot)$ and $Y(\cdot)$. Observe that the $h$ value is obtained
such that $\hat{a}=0.12$.

We consider the goodness-of-fit tests given in Section 5.3 to analyze the appropriate adjustment for these data. Table \ref{TestsApplic2} presents the four goodness-of-fit tests given in Section 5.3 applied to both the original and the twice-differenced time series
for the cardiovascular mortality in Los Angeles County data set. This table considers the statistical tests with their respective $p$-values (in parentheses) for testing hypotheses $H_0$ versus $H_1$, defined in expression \eqref{hypotheses1}. It also presents the corresponding statistical values when one considers the second hypotheses test defined in equation \eqref{hypotheses2}, where $\hat{\alpha}$ is the maximum likelihood estimate value, obtained from the procedure described in Section 5.1.
From the results in Table \ref{TestsApplic2} one observes that all four tests accepted the null hypothesis of
a Gaussian distribution ($\alpha=2$). These all four tests also accepted the null hypothesis
of $\hat{\alpha}=1.9861$ (for $X(\cdot)$) and $\hat \alpha = 1.9998$ (for $Y(\cdot)$). Since
these two maximum likelihood estimate values are very close to $2$,
the results indicate a Gaussian distribution.

\begin{table}[H]\caption{Maximum likelihood estimation results for the original
time series $X(\cdot)$ and the twice-differenced $Y(\cdot)$.}\label{MortalityMLE}
\centering
\begin{tabular}{|c|c|c|}
\hline
\hline
   $\hat{\alpha}$ & $\hat{\sigma}_{\varepsilon}$ &  $\hat{a}$ \Tspace\Bspace\\
\hline
\hline
\multicolumn{ 3}{|c|}{$X(\cdot)$; $h=5.66$}  \Tspace\Bspace\\
\hline
1.9861   & 7.5509   &   0.1209       \Tspace\Bspace\\
\hline
\hline
\multicolumn{ 3}{|c|}{$Y(\cdot)$; $h=18$}  \Tspace\Bspace\\
\hline
1.9998   &    7.4668   & 0.1209   \Tspace\Bspace\\
\hline
\hline
\end{tabular}
\end{table}

\begin{figure}[h!!!]
\centering
\mbox{
\subfigure[]{\includegraphics[width=0.5\textwidth]{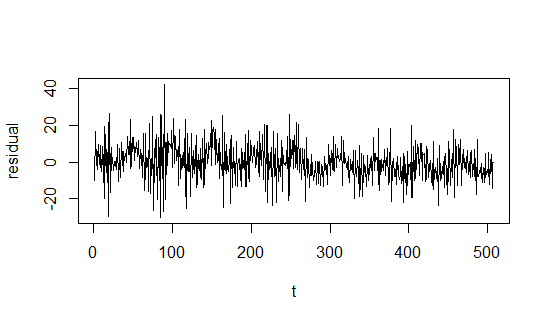}}
\subfigure[]{\includegraphics[width=0.5\textwidth]{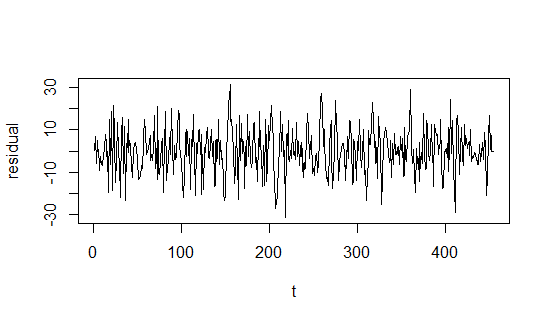}}}
\vspace{-0.3cm}
\caption{Residual analysis for the cardiovascular mortality in Los Angeles County: (a) $X(\cdot)$; (b) $Y(\cdot)$.}
\label{residualsApplic2}
\end{figure}

\begin{table}\caption{Statistics values for the Kolmogorov-Smirnov (KS),
Anderson-Darling (AD), Modified Kolmogorov-Smirnov
(MKS) and McCulloch (MC) tests, applied to both $X(t)$ and $Y(t)$
time series, with their corresponding $p$-values (in parentheses).}\label{TestsApplic2}
\centering
\begin{tabular}{cccc}
\hline
\hline
\mbox{KS} & \mbox{AD} & \mbox{MKS} & \mbox{MC} \\
\hline
\hline
\multicolumn{4}{c}{$X(\cdot)$: $H_0$: $\alpha=2$} \\
\hline
0.1013 (0.217) & 13.7358 (0.259) & 21.1807 (0.392) & 0.13 (0.983) \\
\hline
\multicolumn{4}{c}{$X(\cdot)$: $H_0$: $\hat{\alpha}=1.9861$} \\
\hline
0.1015 (0.379) & 13.8139 (0.495) & 21.2373 (0.354) & 0.105 (0.995)  \\
\hline
\hline
\multicolumn{4}{c}{$Y(\cdot)$: $H_0$: $\alpha=2$} \\
\hline
0.1037 (0.176) & 12.8226 (0.14) & 19.2001 (0.742) & 0.974 (0.071)   \\
\hline
\multicolumn{4}{c}{$Y(\cdot)$: $H_0$: $\hat{\alpha}=1.9998$} \\
\hline
0.1037 (0.178) & 12.8234 (0.123) & 19.2015 (0.744) & 0.968 (0.046)  \\
\hline
\hline
\end{tabular}
\end{table}

\section{Conclusions}\label{conclusionsection}

In this work, we present a novel class of stochastic process that generalizes the
Ornstein-Uhlenbeck processes. This novel class is hereafter called by
\emph{Generalized Ornstein-Uh\-len\-be\-ck Type Process} and it is denoted by GOU type process.
We consider them driven by several classes of noise process such as Brownian motion, symmetric $\alpha$-stable L\'evy motion,
a L\'evy process and a Poisson process. In Sections $2$ and $3$ we presented the
main theoretical results, mainly, we set necessary and sufficient conditions,
under the memory kernel function $\rho(\cdot)$, for the time-stationary and the Markov properties for these processes.

When the GOU type process is driven by a L\'evy noise we proved it is infinitely divisible showing its generating triplet.
Section $4$ we dedicated to the presentation of several examples derived from the GOU type process: for each of them we showed some of their
basic properties as well as some graphical time series realizations. These examples also presented their theoretical
and empirical autocorrelation or normalized codifference functions depending on whether the process has a finite or infinite second moment.

Furthermore, Section $5$ presented two estimation procedures: the classical maximum likelihood and the Bayesian
estimation for the so-called \emph{Cosine process}, a particular process in the class
of GOU type processes presented in Examples $4.3$ and $4.5$: in Section $5.1$ we considered the maximum likelihood estimation
procedure for the \emph{Cosine process} driven both by Brownian motion and symmetric $\alpha$-stable L\'evy noise while in Section $5.2$ we considered the Bayesian estimation procedure for the same \emph{Cosine process} driven by the same two noise processes.
For the Bayesian estimation method, we proposed a power series representation, based on Fox's H-function, to better approximate the density function of a random variable $\alpha$-stable distributed.
In Section 5.3 we present four goodness-of-fit tests (three of them from the EDF class and the McCulloch, 1986 test) for helping to decide which distribution best fits a particular data set.

Finally, in Section $6$ we presented two applications of GOU type model: one is based on the Apple company stock market price data and the other is based on the cardiovascular mortality in Los Angeles County data. The first data set is a time series from January 04, 2010, to November 26, 2014, with an overall total of $481,650$ observations.
We considered four log-return time series: one-minute $R_t^{(1)}$, five-minute $R_t^{(5)}$, ten-minute $R_t^{(10)}$ and fifteen-minute $R_t^{(15)}$
log-returns. The maximum likelihood estimates for the stability index were, respectively, $1.5132$, $1.4661$, $1.4531$, and $1.4304$, and three EDF goodness-of-fit tests showed that
the log-return time series are better fitted by a non-Gaussian $\alpha$-stable process. From the residual values obtained for the adjusted Cosine
model, we concluded that the original log-return time series $R_t^{(1)}$ captures well the features of this data set. This log-return time series has the smallest discretization step size and presented the smallest standard deviation estimate and
the smallest variation of the residual values. The second data set is a weekly time series
given the cardiovascular mortality in Los Angeles County from 1970 to 1979, with an
overall total of $508$ observations. We considered the original and the twice-differenced
time series. The maximum likelihood estimate values for the stability index were, respectively, $1.9861$ and $1.9998$ that are both very close to $2$. Four goodness-of-fit tests accepted the null hypothesis of a Gaussian distribution ($\alpha=2$). These all four tests also accepted the null hypothesis of $\hat{\alpha}=1.9861$ (for the original time series) and $\hat \alpha = 1.9998$ (for the twice-differenced time series). These results
indicate a Gaussian distribution. From the residual values obtained for the adjusted Cosine
model, given in Example $4.3$, we concluded that the twice-differenced time series has the
smallest estimate value for the standard deviation and their residuals behave like a random
noise around the zero mean value.

\subsection*{Acknowledgments}

J. Stein was supported by CNPq-Brazil.
S.R.C. Lopes's research was partially supported by CNPq-Brazil.
A.V. Medino's research was partially supported  by FAPDF through the research grant
No. 193.000.061/2012: ``Infer\^encia em Processos Estoc\'asticos de Alta Variabilidade e de Longa Depend\^encia".


\end{document}